\documentclass{siamart220329}	

\usepackage{array,boldline,makecell} 
\usepackage[utf8]{inputenc}
\usepackage{graphicx,multirow}  
\usepackage{verbatim,float} 
\usepackage{url,todonotes}
\usepackage{amssymb,amsmath,mathtools}
\usepackage{graphicx,amsmath,amsfonts,amssymb,subfigure}
\usepackage{color}

\usepackage{algorithm}
\usepackage{algpseudocode}
\usepackage{textcomp}
\def\RR{\mathbb{R}}

\def\XX{\mathcal{X}}
\def\FF{\mathcal{F}}

\def\Span{\mbox{Span}}

\def\betab{\vspace{-0.05in}\begin{tabbing} 
xxx\=xxxx\=xxxx\=xxxxxx\=xxxxxx\=xxxxxxxxxxxxxxxxxxxxxxxx\=xxxxxxxxx\= \kill} 
\def\entab{\end{tabbing}}

\def\Tr{\text{tr}}

\def\RR{\mathbb{R}}

\def\up#1{^{({#1})}} %
\def\nref#1{(\ref{#1})}
\def\inv{^{-1}}%
\def\half{\frac{1}{2}} 
\def\eps{\epsilon}%

\newcommand{\eq}[1]{\begin{equation}\label{#1}}
\newcommand{\en}{\end{equation}}

\graphicspath{{./FIGS/}}

\title{
  NLTGCR: a class of nonlinear
  acceleration procedures based on Conjugate Residuals}
\author{Huan He\thanks{Work done in Department of Computer Science, Emory University, Atlanta, GA 30322 (\email{hehuannb@gmail.com})}
\and Ziyuan Tang\thanks{Department of Computer Science and Engineering, University of Minnesota, Minneapolis (\email{tang0389@umn.edu}, \email{saad@umn.edu}). The research of Tang and Saad is supported by the NSF award DMS 2208456.}
\and Shifan Zhao\thanks{Department of Mathematics, Emory University, Atlanta, GA 30322 (\email{shifan.zhao@emory.edu}, \email{yxi26@emory.edu}). The research of Zhao and Xi is supported by NSF award  DMS 2208412.}
\and Yousef Saad\footnotemark[2]
\and Yuanzhe Xi\footnotemark[3]}

\headers{Nonlinear Truncated GCR}{H. He, Z. Tang, S. Zhao, Y. Saad, and Y. Xi} 

\begin{document} 

\maketitle 

\begin{abstract}
  This paper develops a new class of nonlinear acceleration algorithms based on
  extending conjugate residual-type procedures from linear to nonlinear
  equations. The main algorithm has strong similarities with
  Anderson acceleration as well as with inexact Newton methods - depending
  on which variant is implemented. We prove theoretically and verify experimentally, on
 a variety of problems from simulation experiments to deep learning applications, that
  our method is a powerful accelerated iterative algorithm. {The code is  available at \url{https://github.com/Data-driven-numerical-methods/Nonlinear-Truncated-Conjugate-Residual}}.
  
\end{abstract}

\begin{keywords} 
  Nonlinear acceleration, Generalized Conjugate Residual, Truncated GCR,
  Anderson acceleration, Newton's method, Deep learning
\end{keywords}

\begin{AMS}
   	65F10, 68W25, 65F08, 90C53
\end{AMS}
\section{Introduction}\label{sec:intro}
There has been a surge of interest in recent years in numerical algorithms
whose goal is to  accelerate  iterative schemes for solving the following problem: 
\eq{eq:Fx} \text { Find } x \in \mathbb{R}^{n} \text { such that } f(x) = 0 ,
\en
where $f$ is a continuously differentiable mapping from $\RR^n $ to $\mathbb{R}^{n}$.
This  problem can itself originate from  unconstrained optimization
where we need to minimize a scalar function:
\eq{eq:Opt} \min_x \phi(x) , \en in
which $\phi: \ \RR^n \to \RR$.  In this situation, we will be interested in a
local minimum which can be found as a zero of the system of equations $f(x) = 0$ where
$f(x) = \nabla \phi(x)$.

The problem \eqref{eq:Fx} can be formulated as a \emph{fixed point} problem, where one seeks to
find  the fixed point of a mapping $g$ from $\mathbb{R}^{n}$ to itself:
\begin{equation}\label{eq:fp}
\text { Find } x \in \mathbb{R}^{n} \text { such that } x = g(x) . 
\end{equation}
This can be achieved by setting, for example, $g(x) = x + \beta f(x)$ for some
nonzero $\beta$.  Given a mapping $g$, the related \emph{fixed point iteration},
i.e., the sequence generated by \eq{eq:fp1} x_{j+1} = g (x_j) \en may converge
to the fixed point of \eqref{eq:fp}  and when $g(x) = x+\beta f(x)$ then this limit
is clearly a solution to the problem \eqref{eq:Fx}.  However, fixed-point
iterations often converge slowly, or may even diverge. As a result
\emph{acceleration methods} are often invoked to improve their convergence or
to establish it.


A number of such acceleration methods have been proposed in the past.
It is important to clarify the terminology and discuss the  distinction between
\emph{acceleration methods} which aim at accelerating the convergence of a fixed point sequence
of the form \eqref{eq:fp1}, and \emph{solvers} which aim at finding solutions to \eqref{eq:Fx}.
Among acceleration techniques are
`extrapolation-type' algorithms such as the Reduced-Rank Extrapolation (RRE)
\cite{doi:10.1137/1029042}, the Minimal-Polynomial Extrapolation (MPE)
\cite{doi:10.1137/0713060}, the Modified MPE (MMPE) \cite{8145128}, and the
vector $\epsilon$-algorithms \cite{doi:10.1137/140957044}. These typically
produce a new sequence from the original one by combining them without invoking
the mapping $g$ in the process.
Another class of methods produce a new sequence
by utilizing both the iterates and the mapping $g$.  Among these, Anderson
Acceleration (AA) \cite{Anderson65} has received enormous attention in recent
years due to  its success in solving a wide range of problems 
\cite{09639,AARL,Shi2019RegularizedAA,d_Aspremont_2021,Sun2021DampedAM,9763953,he2022gdaam, he2020fast, wei2022a, he2021age}.
AA can be seen as an inexpensive alternative to second order methods
such as quasi-Newton type techniques.  It is often used
quite successfully without global convergence strategies such as line search or
trust-region techniques. These advantages made the method popular in
applications ranging from quantum physics, where they were first
developed, to machine learning.
We refer readers to \cite{shanks} for a survey of acceleration methods.

As was stated above, it is clear that one can invoke one of these accelerators
for solving \nref{eq:Fx} by applying it to the fixed point sequence associated
with $g(x) = x + \beta f(x)$.  AA, and its sibling Pulay mixing \cite{pul80,Pulay-DIIS},
were devised
precisely in this way.  Thus, \emph{an accelerator of this type can be viewed as
  a solver}, in the same way that a linear accelerator (e.g., Conjugate
Gradient) combined with some basic iteration, such as Richardson, can be viewed
as a `solver'.  This class of techniques does not include the extrapolation-type
methods discussed above because \emph{they require the computation of $g(x)$ for
  an arbitrary $x$}. In extrapolation methods we have a sequence of vectors, but
we do not have access to the function $g$ for evaluating $g(x)$.

We can also ask the question of whether or not a given solver can be viewed as an
accelerator.  If the solver only requires evaluating $f(x)$ for an arbitrary $x$
then clearly we can apply it to find the root of the equation
$f(x) \equiv x-g(x)= 0$, which requires the computation $g(x)$ given $x$. The
related iterative process can be viewed as an acceleration technique for the
fixed point mapping $g(x)$.
Thus, our definition of an accelerator is broad and it encompasses any method
that aims to speed up a fixed point iteration by requiring only function
evaluations at each step.

A good representative of this class of methods is Anderson Acceleration (AA).  There
are three issues with AA, and similar accelerators, which this article aims to
address. The first is that for optimization problems, AA does not seem to be
amenable to exploiting the symmetry of the Jacobian or Hessian.  If we had to
solve a linear system with AA, the sequence resulting from the algorithm cannot
be written in the form of a short-term recurrence, as is the case with Conjugate
Gradient or Conjugate Residual algorithm for example. For nonlinear optimization
problems where the Hessian is symmetric, this indicates that AA does not take
advantage of symmetry and as such it may become expensive in terms of memory and
computational cost.  This is especially true in a nonconvex stochastic setting,
where a large number of iterates are often needed.  This can be an acute
problem, particularly in machine learning, where we often encounter practical
situations in which the number of parameters is quite large, making it
impractical to use a large number of vectors in AA.  Although recent literature
made efforts to improve the convergence speed of AA
\cite{wei2021stochastic,wei2022a,8682962,09639,03971,he2022efficient}, they did not attempt to reduce memory
cost.  One of the goals of this paper is to propose an acceleration technique
that exploits symmetry or near symmetry to reduce the computational cost.

A second problem with AA, which is perhaps a result of its simplicity, is that
while it shares some similarities with quasi-Newton methods, it does not exploit
standard methods that are common in second order methods to improve \emph{global
  convergence} characteristics. The inclusion of line search or trust-region
methods is necessary if one wishes to solve realistic problems.  A secondary
goal of this paper is to introduce an AA-like method 
that implements global convergence strategies borrowed from inexact Newton and
quasi-Newton methods.

Finally, one of the issues with AA, and other acceleration methods, is that
it uses a crude linear model.  Specifically, acceleration methods typically
rely on two sets of consecutive vector differences, namely the differences
$\Delta f_i \equiv f(x_{i+1}) - f(x_{i}) $ and the associated
$\Delta x_i \equiv x_{i+1} -x_{i} $. The issue is that their approximate
solutions are developed from the relation
$\Delta f_i \approx J(x_i) \Delta x_i$, where $J(x_i)$ is the
Jacobian at $x_i$. This linear approximation is
likely to be rather inaccurate, especially when the iterate is far from its
limit. It is important to develop techniques that will avoid  relying on such rough
approximations.

Nonlinear acceleration methods of the type discussed in this paper appeared
first in the physics literature where they were needed to accelerate very
complex and computationally intensive processes, such as the Self-Consistent
Field (SCF) iteration. The best-known of these methods was discovered by
Anderson~\cite{Anderson65} in 1965. In the early 1980s, Pulay proposed a similar
scheme which he called Direct Inversion on the Iterative Subspace
(DIIS)~\cite{pul80,Pulay-DIIS}.  Both methods were designed specifically for SCF
iterations and it turns out that, although formulated differently, AA and DIIS
are essentially equivalent, see, e.g.,
~\cite{chupin2021convergence,Fattebert10Accel,LinChao13}.
For this reason, the method is often referred to as `Anderson-Pulay mixing',
where mixing in this context refers to the process by which a new fixed point
iterate is \emph{mixed} with previous ones to accelerate the process. 
AA was rediscovered again in a different form in a 2000 paper by Oosterlee and Washio
\cite{OosterleeWashio} 
who applied their technique to accelerate nonlinear multigrid iterations. 

The link between nonlinear acceleration methods such as AA and secant-type
methods was first unraveled by Eyert \cite{eyert:acceleration96} when he 
compared AA with a multi-secant method proposed by Vanderbilt and Louie more than a
decade earlier \cite{vl:energy84}.  The article~\cite{FangSaad07} explored this idea further
and expanded on it by proposing two classes of multi-secant methods.
Thereafter, AA started to be studied and utilized by researchers outside the field of
physics, see, e.g.,
~\cite{degr,WalkerNi2011,TothKelley15,pmlr-v119-mai20a,Bian-alAA21,Pollock19,Pollock20}
among many others. 

The primary contribution of the present paper is to take another look at this
class of methods and develop a technique that is derived by a careful extension
of a linear iterative method to nonlinear systems.  The paper is motivated
primarily by a desire to overcome the three weaknesses of AA mentioned earlier and { accelerate the stochastic optimization algorithms used in deep learning applications. }

\section{Background: Inexact Newton, quasi-Newton, and Anderson acceleration} 
The goal of this section is to clarify key features of the method proposed
in this paper in order to establish links with known methods.
Many of the approaches developed for solving \eqref{eq:Fx} 
are rooted in Newton's method which exploits the local linear model:
\begin{equation}
\label{eq:linmod}
f(x+\Delta x)\approx f(x) + J(x)\Delta x,
\end{equation}
where $J(x)$ is the Jacobian matrix at $x$.

\paragraph{Notation} We will often refer to an evolving set of columns
where the most recent vectors from a sequence are retained. In this situation,
we found it convenient to  use the following convention. For a given $m\ge 0$
we set:
\eq{eq:not}
j_m = \max\{ 0, j-m+1\}
\quad
m_j = \min \{ m, j+1 \} \equiv j - j_m +1 . 
\en

\subsection{Inexact Newton methods} \label{sec:InexNewt}
Newton's method determines  $\Delta x_j = x_{j+1} - x_j$  
at step $j$, to make the right-hand side on \eqref{eq:linmod}
equal to zero when $x = x_j$, which is achieved by solving
the Newton linear system $ J(x_j) \delta + f(x_j) = 0.$ 
Inexact Newton methods, e.g., \cite{Kelley-book,Dembo-al,Brown-Saad}
among many references, compute a sequence of iterates in which the
above Newton systems are solved approximately, often by an
iterative method.  Given an initial guess $x_0$, the iteration 
proceeds as follows:
\begin{align} 
 &      \mbox{Solve}  & J(x_j) \delta_j   &\approx -f(x_j) &  \label {eq:inex1}\\
 &       \mbox{Set}    & x_{j+1}  & = x_j + \delta_j & \label {eq:inex2}
  \end{align}  
  Note that the right-hand side of  the Newton
  system is $-f(x_j)$ and this is also
  the residual for the linear system when $\delta_j = 0$. 
  Therefore, in later sections we will define the residual vector $r_j$  to  be $r_j \equiv -f(x_j)$.

   The technique for solving the local system \nref{eq:inex1} is not specified.
  Suppose that we invoke a Krylov subspace method for solving~\nref{eq:inex1}. If we set $J \equiv J(x_j)$ then the method,
  will usually generate an approximate solution that
  can be written in the form
    \eq{eq:deltaj}
  \delta_j = V_j y_j,
  \en
  where $V_j$ is an orthonormal basis of the Krylov subspace
  \eq{eq:Kry}
  K_j = \Span \{r_j, J r_j, \cdots, J^{m-1} r_j \}.
  \en
 The vector $y_j$ represents the expression of
  the solution in the basis $V_j$. 
For example, if GMRES or, equivalently Generalized Conjugate Residual (GCR) \cite{Eis-Elm-Sch}, is used, then $y_j$ becomes $y_j = (J V_j)^\dagger (-f(x_j))$, {where $\dagger$ denotes the pseudoinverse.} In essence the inverse Jacobian is approximated by the rank $m$ matrix:
\[B_{j, GMRES} = V_j (J V_j)^\dagger. \] 
In inexact Newton methods the
approximations just defined are valid only for the $j$-th step, i.e.,
once the solution is updated, they are discarded and the process will
essentially compute a new Krylov subspace and related approximate
Jacobian at point $x_{j+1}$. This `lack of memory' can be an
impediment to performance. In contrast, quasi-Newton methods will
compute approximate Jacobians or their inverses by a process that is
continuously being updated, using  the most recent iterate
for this update.

  \subsection{Quasi-Newton methods}
\label{sec:broyden}
Standard quasi-Newton methods build a  local approximation 
$J_{j}$ to the Jacobian $J(x_{j}) $ progressively by using previous iterates.  These methods require  the relation  \eqref{eq:linmod} to be satisfied by
the updated $J_{j+1}$ which is built at step $j$. 
This means that the following {\em secant condition}, is imposed:
\begin{equation}
\label{eq:secant1}
J_{j+1}\Delta x_j = \Delta f_j,
\end{equation}
where $\Delta f_j:=f(x_{j+1})-f(x_j)$.
The following
{\em no-change condition} is also imposed:
\begin{equation}
\label{eq:nochange1}
J_{j+1}q =  J_{j}q, \quad \forall q \quad \mbox{such that} \quad 
q^T\Delta x_j=0 . 
\end{equation}
In other words, there should be no new information from $J_j$ to $J_{j+1}$
along any direction $q$ orthogonal to $\Delta x_j$.
Broyden showed that there is a unique  matrix $J_{j+1}$ that satisfies
both  conditions (\ref{eq:secant1}) and (\ref{eq:nochange1}) and it can be
obtained by the update formula: 
\begin{equation} 
J_{j+1} = 
J_j + (\Delta f_{j}-J_j\Delta x_j)\frac{\Delta x_j^T}{\Delta x_j^T\Delta x_j}.
\label{eq:broyden1updatej}
\end{equation}


Broyden's \emph{second method} approximates
the inverse Jacobian directly instead of the Jacobian itself.
If $G_j$  denotes this approximate inverse Jacobian
at the $j$-th iteration, then the secant condition (\ref{eq:secant1})
becomes: 
\begin{equation}
\label{eq:secant2}
G_{j+1}\Delta f_j = \Delta x_j.
\end{equation}
By minimizing $E(G_{j+1})=\|G_{j+1}-G_{j}\|_F^2$
with respect to $G_{j+1}$
subject to (\ref{eq:secant2}),
one  finds this update formula for the inverse Jacobian:
\begin{equation}
\label{eq:broyden2update}
G_{j+1}=G_{j}+(\Delta x_j - G_j 
\Delta f_j)\frac{\Delta f_j^T}{\Delta f_j^T\Delta f_j},
\end{equation}
which is also the only update satisfying both
the secant condition (\ref{eq:secant2})
and the no-change condition for the inverse Jacobian:
\eq{eq:NoCh}
(G_{j+1} -G_{j}) q = 0, \quad \forall \  q \perp \Delta f_j.
\en

We will revisit secant-type methods again when we discuss AA in the next section. AA can be
viewed from the angle of \emph{multi-secant} methods, i.e., block forms of the secant
methods just discussed, in which we impose a secant condition on a whole set of vectors
$\Delta x_i, \Delta f_i$ at the same time.

\subsection{General nonlinear acceleration and Anderson's method} 
\label{sec:AA1} 
Acceleration methods, such as AA, take a different viewpoint altogether. Their
goal is to accelerate a given fixed point iteration of the form
\eqref{eq:fp1}.
Thus, AA starts with an
initial $x_0$ and sets  $x_1=g(x_0)=x_0+\beta f_0$, where
$\beta >0$ is a parameter.  At step $j>m$ we define
$X_j=[x_{j-m},\ldots, x_{j-1}],$   and
$ F_j=[ f_{j-m}, \ldots, f_{j-1}]$ along with the differences:
\begin{equation}
\label{eq:dfdx}
\mathcal{X}_j=[\Delta x_{j-m}\;\cdots\;\Delta x_{j-1}]\in \RR^{n\times m},
\qquad
\mathcal{F}_j=[\Delta f_{j-m}\;\cdots\;\Delta f_{j-1}]\in \RR^{n\times m}.
\end{equation}
We then  define the next AA iterate as follows: 
\begin{align} 
x_{j+1} &=  x_j+\beta f_j -(\XX_j+\beta \FF_j)\ \theta\up{j}  \quad \mbox{where:}  \label{eq:AA}  \\
\theta\up{j} &= \text{argmin}_{\theta \in \mathbb R^{m}}\| f_j - \FF_j \theta \|_2 .  
 \label{eq:thetaj}
\end{align}
Note that $x_{j+1}$ can be expressed with the  help of intermediate vectors:
\eq{eq:AA1} 
\bar x_j = x_j-\XX_j \ \theta\up{j} ,\quad 
\bar f_j = f_j-\FF_j \ \theta\up{j} , \quad 
x_{j+1} =\bar x_j + \beta \bar f_j .
\en

There is an underlying quasi-Newton second order method aspect to the procedure.
In Broyden-type methods,
Newton's iteration: $x_{j+1}  = x_j - J(x_j) \inv f_j$ is replaced with
$     x_{j+1} = x_j - G_j f_j $  
  where $G_j $  approximates the inverse of the Jacobian $J(x_j)$ at $x_j$
  by the  update formula $G_{j+1} = G_j + (\Delta x_j - G_j \Delta f_j) v_j^T $
  in which $v_j$ is defined in different ways see \cite{FangSaad07} for details. 
AA belongs to the class of \emph{multi-secant methods}.
Indeed, the approximation \nref{eq:AA} can be written as: 
\eq{eq:AndQN} 
x_{j+1} = x_j - [- \beta I  + (\XX_j+\beta\FF_j)(\FF_j^T\FF_j)^{-1}\FF_j^T ] f_j\equiv  x_j - G_j f_j .
\en 
Thus,  $G_{j} $ can be seen as an update to
the (approximate) inverse  Jacobian $ G_{j-m} = -\beta I$  by the formula:
\eq{eq:msec}
G_{j}  = G_{j-m} + (\XX_j  - G_{j-m} \FF_j)(\FF_j^T\FF_j)^{-1}\FF_j^T. 
\en 
It can be shown that the approximate inverse
Jacobian $G_j$ is the result of 
minimizing $\| G_j + \beta I  \|_F$ under the \emph{multi-secant condition}
of type II \footnote{Type I Broyden conditions involve
  approximations to  the Jacobian,
 while type II conditions  deal with the inverse Jacobian.}
\eq{eq:mscond}
G_j \FF_j = \XX_j.
\en
This link between AA and Broyden multi-secant type updates
was first unraveled by Eyert~\cite{eyert:acceleration96} and expanded upon in
\cite{FangSaad07}. Thus, the method is in essence what we might call a
`block version' of Broyden's second update method, whereby  a rank $m$,
instead of rank $1$, update is applied at each step.

\subsection{The issue of symmetry} 
  Consider again the specific case where the nonlinear function
$f(x)$ is the gradient of some
scalar function $\phi(x)$ to be minimized, i.e., $f(x) = \nabla \phi(x) $.
In this situation the Jacobian of $f$ becomes $\nabla^2 \phi$  the Hessian of $\phi$,
  and therefore it is symmetric.
 Approximate Jacobians that are implicitly or explicitly extracted in the
  algorithm, will be  symmetric or nearly symmetric. 
  Therefore this raises the
  possibility of developing accelerators that take advantage of
  symmetry or near-symmetry. One way to achieve this is to extend
  \emph{linear solvers} that take advantage of symmetry to the nonlinear context.  This was one of the primary initial motivations
  for this work.

  We saw earlier that AA is a multi-secant version of a
  Broyden type II method where the approximate inverse Jacobian is updated by formula
  \nref{eq:msec}.  An obvious observation here is that the symmetry of the Jacobian can not be exploited in any way in this formula.  This has been considered
  in the literature (very) recently, see for example,
  \cite{Boutet2020,scieur21a,Boutet2021}.  In a 1983 report, Schnabel
  \cite{Schnabel83} showed that the matrix $G_j$ obtained by a multi-secant method
  that satisfies the secant condition \nref{eq:mscond} is symmetric iff the
  matrix $\XX_j^T \FF_j$ is symmetric. It is possible to explicitly force
  symmetry by employing generalizations of the symmetric versions of
  Broyden-type methods. Thus, the authors of \cite{Boutet2020,Boutet2021}
  recently developed a multi-secant version of the Powell Symmetric Broyden (PSB)
  update due to Powell \cite{Powell70} while the article \cite{scieur21a}
  proposed a symmetric multi-secant method based on the popular
  Broyden-Fletcher-Goldfarb-Shanno (BFGS) approach as well as the
  Davidon-Fletcher-Powell (DFP) update. However, there are a number of  issues with the  symmetric versions of multi-secant updates
  some of which are  discussed in~\cite{scieur21a}.

\section{Nonlinear Truncated Generalized Conjugate Residual (nlTGCR) algorithm}\label{sec:linear} 
For understanding the conjugate residual-based methods in the nonlinear
case, it is important to first provide some background for
linear systems. A large
class of Krylov subspace methods for solving nonsymmetric linear systems have
been developed in the past four decades. The reader is referred to the recent
volume by Meurant and Tebbens \cite{Meurant-Tebbens} which contains a rather
exhaustive and detailed coverage of these methods.
The main aim of the techniques proposed in
this article is to  adapt the residual-minimizing
subclass of Krylov methods for linear systems to the nonlinear case.
The guiding principle in this adaption is that we 
would like it to also approximately minimize the
nonlinear residuals. This is in contrast with inexact Newton methods
where the goal is to roughly solve the linear systems that arise in
Newton's method as a way to provide a good search direction.

\subsection{The linear case:  Generalized Conjugate Residual (GCR) Algorithm}\label{sec:gcr} 
We first consider solving a linear system of the form: \eq{eq:Ax=b} A x = b.
\en A number of iterative methods developed in the 1980s aimed at minimizing the
norm of the residual $r = b -A x$ of a new iterate that lies in a Krylov
subspace, see \cite{Meurant-Tebbens} for a detailed account.  Among these, we
focus on the Generalized Conjugate Residual (GCR) algorithm \cite{Eis-Elm-Sch}
for solving \nref{eq:Ax=b} which is sketched in Algorithm \ref{alg:gcr}.

    \begin{algorithm}[htb]
    \centering
    \caption{GCR}\label{alg:gcr}
    \begin{algorithmic}[1]
  \State \textbf{Input}: Matrix $A$, RHS $b$,
  initial  $x_0$. \\
  Set $p_0 := r_0 \equiv b-Ax_0$.
\For{$j=0,1,2,\cdots,$ until convergence} 
\State  $\alpha_j :=  \langle r_j, A p_j \rangle / \langle A p_j, A p_j\rangle $
\State $x_{j+1} := x_j + \alpha_j p_j$
\State $r_{j+1} := r_j - \alpha_j A p_j$
\State $p_{j+1} := r_{j+1} - \sum_{i=0}^j \beta_{ij} p_i $ \quad where \quad
$\beta_{ij} :=  \langle A r_{j+1} , Ap_i\rangle / \langle A p_i, Ap_i\rangle $
\EndFor
\end{algorithmic}
\end{algorithm}
The main point of the algorithm is to build a sequence of search
directions $p_i$, $i=0:j$ at step $j$ so that the vectors $A p_i$ are
orthogonal. This is done in Line~7.  With this we know that the
iterate as defined by Lines 4-5 is optimal in the sense that it yields
the smallest residual norm among solution vectors selected from the
Krylov subspace $x_0 + \mbox{\Span}\{p_0, \ldots, p_k\}$ -- see \cite[pp 195-196]{Saad-Book2}.  The GCR
algorithm is mathematically equivalent
\footnote{Here equivalent is meant in the sense that if exact arithmetic is
  used and if the compared algorithms both succeed in  producing the $j$-th
  iterate from the same initial $x_0$, then the two iterates are equal.} 
to GMRES~\cite{Saad-book3}, and
to some of the forms of other methods developed earlier, e.g.,
ORTHOMIN~\cite{Vinsome76}, ORTHODIR~\cite{Jea-Young}, and Axelsson's
CGLS method \cite{Axelsson80}.  

Next, we will discuss a truncated version of GCR
in which the $Ap_i$'s are only orthogonal to the $m$ previous
ones instead of all of them.
This algorithm was first introduced by Vinsome as early as 1976
and was named `ORTHOMIN' \cite{Vinsome76}.
We will just refer to it as the truncated version of GCR or TGCR(m).
In a practical implementation, we need to keep a set of $m$ vectors 
at step $j$ for the $p_i$'s and another set for the vectors $v_i = A p_i$.
In addition, we replace the classical Gram-Schmidt of Line~7 of Algorithm~\ref{alg:gcr}
by the modified form of Gram-Schmidt: 
 the vector $A r_{j+1} $ initially set to a vector $v$ which is
 orthogonalized against the successive $A p_i$'s. Thus,
 Line~7 of Algorithm~\ref{alg:gcr} becomes: 

 \medskip 
 \begin{tabbing}
   xxxxx\=xxxx\=xxxxxxxxxxxxxxxxx\=xxxxxxxx\=xxxx\=xxxx\kill
 7a.  \> $p := r_{j+1}$; $ v := A p $;  and  $j_m := \max(0,j-m+1)$  \\
 7b. \> \textbf{for} \ $i=j_m:j$ \\
 7c. \> \>  $\beta_{ij} :=  \langle v, Ap_i \rangle $\\
 7d. \> \> $p := p - \beta_{ij} p_i$; \>  $v := v - \beta_{ij} v_i$; \\
 7e. \> \textbf{end for}\\
7f. \>  $p_{j+1} :=p /\| v\|$ ; \>\>   $v_{j+1} :=v/\|v\|$ ; 
\end{tabbing}

\medskip
We refer to the Algorithm obtained from Algorithm~\ref{alg:gcr} where
Line~7 is replaced by Lines (7a--7f) above as the Truncated GCR (TGCR(m)) algorithm.
This algorithm is a slight variant of the original ORTHOMIN introduced 
in \cite{Vinsome76} and analyzed in  \cite{Eis-Elm-Sch}.
It  produces a system of vectors
$V_{j+1} = [v_{j_m}, v_{j_{m}+1},\ldots,v_j, v_{j+1}]$ that is orthonormal.
When $j \ge m$, $V_{j+1}$ consists of a `window' of $m+1$ vectors. TGCR(m)
approximately  minimizes the quadratic form
$\phi_q(x) \equiv \half \| b - A x\|_2^2 $ in a certain Krylov subspace.  With
$m = \infty $ we obtain the non-restarted GCR method - which is equivalent to
the non-restarted GMRES \cite{Saad-book3}.

A few properties of the vectors generated in TGCR(m) have been analyzed in
~\cite[Th. 4.1]{Eis-Elm-Sch} which  also discussed the convergence of the algorithm
when $A$  is positive definite, i.e., when $A+A^T$ is symmetric positive definite.
When $A$ is symmetric a number of simplifications
take place in Algorithm ~\ref{alg:gcr}. 
In this situation,   all the $\beta_{ij}$'s except 
$\beta_{jj} $ vanish.  
The resulting simplified algorithm yields the standard Conjugate Residual (CR)  algorithm
which dates back to Stiefel \cite{Stiefel-CR}, see
\cite[Section 6.8]{Saad-book3} for details. {See Appendix \ref{sec:Appendix} for a unified presentation of a number of theoretical results of GCR.}

  \subsection{The nonlinear extension: nonlinear TGCR (nlTGCR) algorithm}\label{sec:nonlin}
  We now return to the nonlinear problem and ask the question:  how can we
  generalize the algorithms for linear systems of Section~\ref{sec:gcr} for solving nonlinear
  equations? We should begin by stating what are the desired features of this
  extension.
  First, we would like the algorithms to  fall back to TGCR when the problem is linear.
  Second, we would like a method that can be adapted in such a way as to yield
  the inexact Newton viewpoint when desired or a multi-secant (AA-like) approach when
  desired. Third, we would like a method that exploits a more accurate linear model than
  either Newton or a quasi-Newton approach - possibly at the cost of  a few extra function
  evaluations. Finally, we would like the algorithm to be easily adaptable to the very
  common context in which the function $f$ is `fuzzy' as is the case when dealing with
  stochastic methods.

  In our model, we assume that at step $j$ we have a set 
  of (at most) $m$ current `search' directions
  $\{ p_i \} $ for $j_m \le i \le j$  gathered as columns of a matrix $P_j$,
  where we recall the notation $j_m \equiv \max\{0,j-m+1 \}$.
  {Along with $p_i$'s, we also have a matrix denoted by $V_j$, such that}
    \eq{eq:PjVj}
    P_j = [ p_{j_m}, p_{j_m+1},\cdots, p_j],
    \qquad 
    V_j = [ v_{j_m}, v_{j_m+1},\cdots, v_j].
    \en
    Note that this pair of matrices plays the same role as the pair
    $\XX_j, \ \FF_j$ defined in \eqref{eq:dfdx} for Anderson acceleration.
  {We then  write the linear model used locally as 
   \eq{eq:nltgcrmodle} f(x_j + P_j y) \approx f(x_j) + V_j y. \en }
  While this is again somewhat similar to what was done for Anderson
  acceleration, we note a very important distinction that may have
  a significant effect on performance: \emph{In Anderson acceleration
    the linear model is simply based on the relation
    $f(x_j-\XX_j \theta) \approx f(x_j) -\FF_j\theta$ whereas {nlTGCR evaluates explicitly  the action of $J(x_i)$ on some vector.} This evaluation can be quite accurate if desired.}
  In contrast, the relation
  $f(x_j-\XX_j \theta) \approx f(x_j) -\FF_j\theta$ can be quite rough,
  especially at the beginning of the process where the vectors $\Delta x_j$ and
 $\Delta f_j$ are usually not small.    
 
 { The projection method will minimize the norm
  $\| f(x_j) + V_j y \|_2 $.} This is achieved
  by determining $y$ in such a way that
  \eq{eq:opt} f(x_j) + V_j y \perp \Span \{V_j\} \
  \rightarrow \ (V_j)^T[f(x_j) + V_j y] = 0 \ \rightarrow \ y = 
  V_j^T r_j
  \en
  where it is assumed the $v_i$'s are orthonormal.

Instead  of (\ref{eq:inex1}--\ref{eq:inex2}) of the inexact Newton approach
we now have an iteration of the form
    \begin{align} 
 &      \mbox{Find}  & y_j &= \mbox{argmin}_y \| f(x_j) + V_j y \|_2 &  \label {eq:inex3}\\
 &       \mbox{Set}  & x_{j+1}  & = x_j + P_j y_j & \label {eq:inex4}
  \end{align}  
  A major distinction between this approach and the standard inexact
  Newton is that the latter exploits approximations to the Jacobian
  around one point in order to build the next iterate.  The iterates
  generated by the iterative process can be viewed as intermediate
  points but they rely on a Jacobian $J(x_0)$ calculated  at the initial
  approximation $x_0$.
  We will revisit the inexact Newton viewpoint in a  later   section.

   The idea of the nonlinear version of the truncated GCR method is to
   exploit the directions that are produced by the TGCR algorithm.
   Note that there is a decoupling between the \emph{update}
   from the current iterate (Lines~(4-5) of GCR/TGCR) to a new one  and the
   \emph{construction of the $p_i$'s} in TGCR (Lines (7a-7f) of TGCR).  In
   essence, the first part just builds a new approximation given a new
   `search' subspace - while the second adds a new direction to this
   evolving subspace.
   This distinction will help us generalize our approach to cases where
   the objective function or the Jacobian varies as the iteration proceeds.
   

    \begin{algorithm}[htb]
    \centering
    \caption{nlTGCR(m)}\label{alg:nltgcr}
    \begin{algorithmic}[1]
  \State \textbf{Input}: $f(x)$,   initial  $x_0$. \\
  Set $r_0 := - f(x_0)$.
  \State Compute  $v := J(x_0)  r_0$; \ (Use Frechet)
  \State $v_0 := v/ \| v \| $, $p_0 := r_0/ \| v \|_2 $;
  \For{$j=0,1,2,\cdots,$ } 
\State  $y_j := V_j^T r_j $ 
\State $x_{j+1} := x_j + P_j y_j$
\State $r_{j+1} := -f(x_{j+1}) $ 
\State Set: $p := r_{j+1}$; and $j_m := \max(0,j-m+1)$
\State Compute $v =  J(x_{j+1})  p$ (Use Frechet)
\For{$i=j_m:j$}
\State $\beta_{ij} :=  \langle v, v_i \rangle $
\State $p := p - \beta_{ij} p_i$
\State $v := v - \beta_{ij} v_i$
\EndFor
\State $p_{j+1} :=p /\| v\|_2$ ; \qquad  $v_{j+1} :=v/\|v\|_2$ ; 
\EndFor
\end{algorithmic}
\end{algorithm}

We now derive our general algorithm from which a few variants will follow.
The algorithm is an extension of the TGCR(m) algorithm discussed above -- with a few
needed changes that reflect the nonlinearity of the problem.
The first change is that any residual is now to be replaced by the negative of
$f(x)$ so $r_0$  becomes   $r_0 = -  f(x_0)$ and Line~6 of GCR/TGCR(m) must be replaced by
  $r_{j+1}  = -  f(x_{j+1})$. In  addition, the matrix $A$ in   the 
  products  $A  r_0  $  and  $A  p$ invoked  in    Line~2  and  Line~7a
  respectively, is the   Jacobian of  $f$ at the most recent 
  iterate.  The most
  important  change is  in Lines~(4-5) of Algorithm \ref{alg:gcr} where  according  to the  above
  discussion the scalar $\alpha_j$   is to be replaced
  by the 
  vector $y_j$ that  minimized $\| f(x_j) + V_j y \|_2$ over $y$.
  The reason for this was explained above.  

The resulting nlTGCR(m) algorithm is shown in Algorithm~\ref{alg:nltgcr}.  
It requires two function evaluations per step: one in Line 8 and the
other in Line~10.  {Alternatively, when constructing the Jacobian
is inexpensive,} one can compute $J p$ in Line~10 as a matrix-vector
product. Clearly, the system $[v_{j_m}, v_{j_m +1}, \cdots,  v_{j+1}] $ is orthonormal. 

\section{Theoretical results}
This section discusses connections of the nlTGCR(m) algorithm with inexact and quasi-Newton methods
and analyzes  its convergence.
\subsection{Linearized update version and the connections to inexact Newton methods} 
\label{sec:LinUpd} 
We first consider a variant of the algorithm which we call
the ``linearized update version''. We will show that this version is equivalent to inexact Newton methods in which the system is approximately solved with TGCR(m). Two changes are made to Algorithm \ref{alg:nltgcr} to obtain this linearized update version. 
First,  in Line~8 we update the
residual by using the linear model, namely, we replace Line~8 by: 
 
\bigskip
\betab
\> 8a:  \> \>  $ r_{j+1} = r_j - V_j y_j$
\entab
\smallskip


The second change is that the Jacobian is not updated in Line~10, i.e.,
$J(x_{j+1})$ in Line 10 is kept constant
and equal to $J(x_0)$. In other words
$v$ is computed as
\bigskip
\betab
\> 10a:  \> \>   $ v = J(x_0) p $ 
\entab
\smallskip 


The algorithm is stopped  when the residual norm $r_{j+1} $
is small enough or the number of steps is exceeded.
In addition, \emph{we consider the algorithm merely as a means of providing
  a search direction} as is often done with inexact Newton methods.
In other words the direction $d = x_{last} - x_0$ is provided
to another function that will use it in an iterative procedure
that includes a line search technique at $x_0$. 

It is easy to see that one iteration defined as a sweep using $j$ substeps of
this linear update version is nothing but \emph{an inexact Newton method in
  which the system \nref{eq:inex1} (with $j=0$) is approximately solved with the
  TGCR(m) algorithm.}  Indeed, in this situation the two main loops,  (Lines 3--8 of
Algorithm GCR/TGCR(m)  and Lines 5-17 of Algorithm~\ref{alg:nltgcr})  are
identical.  Lines 5-17 of Algorithm~\ref{alg:nltgcr} perform $k$ steps of the
TGCR(m) algorithm to solve the linear systems $f(x_0) + J(x_0) P y = 0$.  Within
Algorithm~\ref{alg:nltgcr} the update is written in a progressive form as
$x_{j+1} = x_j + \alpha_j p_j$, note that \emph{the right-hand side does not
  change during the algorithm and is equal to $r_0 = - f(x_0)$.  } In effect the
last iterate, $x_k$, is updated from $x_0$ by adding a vector from the Krylov
subspace or equivalently the span of $P_k$.  As a consequence of this
observation, it turns out that $y_j$ has only one nonzero component, namely the
last one.  Indeed, from ~\cite[Th. 4.1]{Eis-Elm-Sch}, we see that if we replace
$A$ by the Jacobian $J$ at $x_0$ then:
\[
  \langle r_{j+1}, J p_i\rangle = 0  \ \mbox{for} \ i=(j-1)_m, \cdots, j. 
\]

It was shown in \cite{Brown-Saad} that under mild conditions, the update to the
iterate is a descent direction. In addition, the article describes `global
convergence strategies' based on line search and trust-region techniques.  If we
do not apply a global convergence strategy, then the algorithm may have difficulties
converging.

Inexact Newton methods \cite{Dembo-al} 
are often implemented with residual reduction stopping criteria of the form
\eq{eq:forcing}
\| f(x_j) + J(x_j) \delta_j \|_2 \le \eta_j \| f(x_j) \|_2
\en
where  $\eta_j \in [0,1) $ is called the forcing term.
This only means that the iterative procedure that is applied when
approximately solving the linear system \nref{eq:inex1} exits when
the relative residual norm falls below $\eta_j$.
A number of articles established convergence conditions under conditions
based on this framework, see, e.g.,
\cite{Dembo-al,Brown-Saad,Brown-Saad2,Eisenstat-Walker94} among others. 

Probably the most significant disadvantage of inexact Newton methods
is that a large 
number of function evaluations may be needed to build the Krylov
subspace in order to obtain a single iterate, i.e., the  next  iterate. 
After this iterate is computed, all the information gathered at
this step, e.g., $P_k, V_k$, is discarded. This is to be contrasted with
quasi-Newton techniques where the most recent function evaluation
contributes to building an updated approximate Jacobian.

\subsection{Nonlinear update version}
The ``linearized update version" of Algorithm~\ref{alg:nltgcr} discussed above uses
a simple linear model to update the residual $r_j$ in which the Jacobian is not updated.
Next, we consider the ``nonlinear update version" as described in
Algorithm~\ref{alg:nltgcr}. 

Assume a sweep using
$j$ substeps of Algorithm~\ref{alg:nltgcr} has been carried out. Our next result will invoke the  linear residual
\eq{eq:LinRes} \tilde r_{j+1} = r_j - V_j y_j , \en
as well as the deviation between the actual residual $r_{j+1}$ and its
linear version $\tilde r_{j+1}$ at the $(j+1)$-th iteration: 
\eq{eq:ResDev} z_{j+1} = \tilde r_{j+1} - r_{j+1} . \en

We first analyze the magnitude of $z_{j+1}$. {Define
\eq{eq:wi}
w_i =J(x_j)p_i -v_i
\ \mbox{for} \ i=j_m, \cdots, j ;  \quad \mbox{and}\quad 
W_j = [w_{j_m}, \cdots, w_j ].
\en
so that, 
\eq{eq:Jxjpi}
  J(x_j) p_i = v_i + w_i .
  \en
Note that the algorithm suggests that $J(x_i) p_i \approx v_i$ but that
equality is not satisfied except in the linear case
\footnote{We thank Eric de Sturler for catching an incorrect statement in an earlier
  version of this paper related to this observation.}}.
  We also define:
\eq{eq:2ndOrd}
s_j  = f(x_{j+1}) - f(x_j) -  J(x_j) (x_{j+1} - x_j) . 
\en
Recall  from the Taylor series expansion $s_j$ is a second order term relative
to $\| x_{j+1}-x_j \|_2$. 
Then we derive the following bound on the norm of $z_{j+1}$ in the next proposition.

\begin{proposition}\label{prop:2}
  The difference $\tilde r_{j+1} - r_{j+1} $ satisfies the relation:
\eq{eq:ResDiff}
\tilde r_{j+1} - r_{j+1}  =   W_j y_j  + s_j
=   W_j V_j^T r_j  + s_j ,
\en
and therefore:
\eq{eq:ResDiffN}
\| \tilde r_{j+1} - r_{j+1} \|_2   \le \| W_j \|_2  \ \| r_j \|_2 + \| s_j \|_2 . 
\en

\end{proposition}

{
\begin{proof} 
We rewrite the difference as:
\begin{align}
  \tilde r_{j+1} -r_{j+1}
  & = (-f(x_j) - V_j y_j ) + f(x_{j+1})  \nonumber\\
  & = (f(x_{j+1}) - f(x_j))  - V_j y_j  \nonumber\\
  & = J(x_j) P_j y_j + s_j - V_j y_j \label{eq:prop2pf}
\end{align}       
Letting $ y_j = [ \eta_{j_m}, \eta_{j_m+1}, \cdots ]$ we  get
from  \nref{eq:Jxjpi} 
\begin{align}
  J(x_j) P_j y_j
  &= \sum  \eta_i J(x_j) p_i  =   \sum  \eta_i [v_i + w_i]  \nonumber \\
  &= V_j y_j  + W_j y_j . \label{eq:prop2pf2} 
  \end{align}
The proof follows by combining \eqref{eq:prop2pf} and \eqref{eq:prop2pf2} \end{proof}}
 
The Proposition \ref{prop:2} shows that $z_{j+1}$ is a quantity of the second order:
when the process is nearing convergence, 
$\| W_j \|_2 \| r_j \|_2 $ is the product of two first order terms
while $s_j$ is a second
order term according to its  definition \nref{eq:2ndOrd}.
 Furthermore, we can prove the following properties of Algorithm~\ref{alg:nltgcr}.
\begin{proposition}\label{prop:nltgcr} 
    The following properties are satisfied by the vectors produced by
    Algorithm~\ref{alg:nltgcr}:

  \begin{enumerate}

  \item 
    $ ( \tilde r_{j+1}, v_i) = 0  \quad \mbox{for} \quad j_m \le i \le j$,
    i.e., $V_{j}^T \tilde r_{j+1} = 0$;

  \item 
  {$ \| \tilde r_{j+1} \|_2     = \min_y \| -f(x_j) + V_j y\|_2 $;}

  \item $\langle v_{j+1} , \tilde r_{j+1}\rangle
    = \left \langle v_{j+1} \ , r_j \right \rangle $;

  \item $y_j=V_{j}^T r_{j} = \langle v_{j}, \tilde r_{j} \rangle e_{m_j} - V_{j}^T z_{j}  $
    where $ e_{m_j} = [0, 0, \cdots, 1]^T \in \ \RR^{m_j}$.
    

\end{enumerate}

\medskip 
\end{proposition}

\begin{proof} 
  Properties (1) and (2) follow from the definition of the algorithm.
  For Property (3)  first observe  that
  $V_j = [v_{j_m}, v_{j_m+1},\cdots, v_{j-1}, v_j]$
  and that $v_{j+1}$ is made orthogonal against  the $m_j$ vectors    $v_{j_m}$ to $ v_j$,
  so $ v_{j+1}^T  V_j  = 0 $ and
  \[
  v_{j+1}^T  \tilde r_{j+1}   =  v_{j+1}^T [r_j - V_j y_j]
  \ = \   v_{j+1}^T r_j . 
  \]

  It is convenient to prove  Property (4)  for the index $j+1$ instead of $j$. We write:
  $     V_{j+1}^T r_{j+1} =  V_{j+1}^T [\tilde r_{j+1} - z_{j+1} ] $ 
where $z_{j+1} $ was defined in Equation \eqref{eq:ResDev}.
  Recalling Properties (1) and (3), we have that
  $\langle \tilde r_{j+1}, v_i \rangle = 0 $ for $j_m \le i \le j$ - so there is only one nonzero term
  in the product $V_{j+1}^T \tilde r_{j+1}$, namely  $v_{j+1}^T \tilde r_{j+1}$.
This gives the results after adjusting for the change of index.
\end{proof}   

Property (4) in Proposition \ref{prop:nltgcr} indicates that when $z_j$ is small, as when the model is \emph{close to being linear} or when it is nearing convergence,
then  $y_{j} $ will have small components everywhere except for the last component. 
As a result it is also possible to consider a slight variant of the algorithm in which  $y_j$ is truncated so as to contain only its last entry.  We will discuss this variant in detail in Section \ref{subsec:linesearch}.

{
\paragraph{Adaptive update}
The nonlinear update version of nlTGCR generally exhibits greater robustness compared to the linearized update version, particularly during the initial phases. In order to leverage the advantages of reduced function evaluations offered by the linearized update version, we introduce an \emph{adaptive update version}. As indicated by property (4) in Proposition \ref{prop:nltgcr}, we employ a mechanism for checking residuals that controls the transition between the two update versions.
Let $r^{nl}_j$ and $r^{lin}_j$ represent the nonlinear and linear residuals at iteration $j$, respectively:
\begin{equation}
  r^{nl}_{j+1} = -f(x_{j+1}), \ \ 
  r^{lin}_{j+1} = r^{nl}_j - V_j y_j.
\end{equation}
We define the `cosine distance' between the nonlinear and linear residuals as follows: 
\begin{equation}
  \theta_j := 1 - \frac{(r^{nl}_j)^Tr^{lin}_j}{\|r^{nl}_j\|_2\cdot\|r^{lin}_j\|_2}.
  \label{eq:switch}
\end{equation}
We choose a threshold $\theta$ which is $0.01$ in this paper. Then, the linearized update version is engaged when $\theta_j < \theta$, while the nonlinear update version is active when $\theta_j \geq \theta$. Note that, when switching from the linearized version to the nonlinear version, the previous vectors $P_j$ and $V_j$ are reset, initiating a new start for nlTGCR with a beginning point at $x_j$. This mechanism effectively reduces the overall number of function evaluations without compromising convergence or accuracy.
}

 \subsection{Connections to quasi-Newton and multi-secant methods} 
In this section, we show that nlTGCR(m) can be viewed from the alternative angle of
a quasi-Newton/multi-secant approach. In this viewpoint,  the inverse of the Jacobian
is approximated progressively. Because it is the inverse Jacobian that is
approximated, the method is akin to Broyden's second update method.

First  observe that in nlTGCR  the update at step $j$ takes the form:

\[x_{j+1} = x_j + P_j V_j^T r_j = x_j + P_j V_j^T (-f(x_j)) . \]
Thus, we are in effect using a secant-type method.
The approximate inverse Jacobian at step
$j$,  denoted by  $G_{j+1} $ for consistency with common
notation is: 
\eq{eq:Gj}
G_{j+1} = P_j V_j^T .
\en 
{
If we apply this to $v_i$ we get
$
  G_{j+1} v_i = P_j V_j^T v_i =  p_i \quad \mbox{for} \quad j_m \le i \le j.
$}  It therefore satisfies 
the \emph{secant} equation
\eq{eq:msecant1} 
  G_{j+1} v_i = p_i \quad \mbox{for} \quad j_m \le i \le j,
  \en
which is a version of \nref{eq:secant2}   used in    Broyden's second update method.
Here,  $p_{i}$ plays the role of $\Delta x_j$ and $v_i$ plays the role
  of $\Delta f_j$.
In addition, the update $G_{j+1}$ satisfies the `no-change' condition:
\eq{eq:msecant2} 
  G_{j+1} q = 0 \quad \forall q \perp v_i  \quad \mbox{for} \quad j_m \le i \le j.
  \en
  The usual no-change condition for secant methods is of the form
  $(G_{j+1}-G_{j-m}) q = 0 $ for $q \perp \Delta f_j$ which in our case
  would be
  $(G_{j+1}-G_{j-m}) q = 0 $ for $q \perp \ v_i  \quad \mbox{for} \quad j_m \le i \le j$.
  One can  therefore consider that we are updating $G_{j-m} \equiv 0$.

  Thus, consider the optimization problem
  \eq{eq:optB}
  \min \{  \| G \|_F  \ \text{subject to: }\quad G V_j = P_j \}
  \en
  which will yield the matrix of the smallest F-norm satisfying the condition \nref{eq:msecant1}.
  Not surprisingly this matrix is just $G_{j+1}$.
  
  \begin{proposition}
  \label{prop:Boptim}
    The unique minimizer of Problem~\eqref{eq:optB} is the matrix $G_{j+1} $ given by
    \nref{eq:Gj}.
    \end{proposition}
    \begin{proof}
      The index $j$ is dropped from this proof. Exploiting orthogonal projectors
      we write $G$ as follows: $G = G V V^T + G (I - V V^T) = P V^T + G (I - V V^T)$ and observe that
      \begin{align*}
        \| G \|_F^2  &= \Tr \left([P V^T + G (I - V V^T) ] [ V P^T + (I - V V^T) G^T ]\right) \\
        &= \Tr \left(P P^T\right)  + \Tr \left(G (I - V V^T)(I-VV^T)G^T \right) \\
                     & = \| P \|_F^2  + \| G (I - V V^T)\|_F^2 
      \end{align*}
      The right-hand side is minimized when $  G (I - V V^T) = 0$ which means  when
      $G = GV V^T$. Recalling the  constraint $G V = P$ yields the desired result.
      \end{proof}

It is also possible to find a link between the method proposed herein
and Anderson acceleration reviewed in Section~\ref{sec:AA1},
by unraveling  a relation with multi-secant methods.
 Based on Proposition \ref{prop:Boptim}, we know that
\eq{eq:msecant1M}
 G_{j+1} V_j = P_j   .   
  \en  
  This is similar to the multi-secant condition
  $ G_j \mathcal{F}_j ={\cal X}_j $ of Equation
  \nref{eq:mscond} -- see also
  Equation (13) of \cite{FangSaad07} where   $  {\cal F}_j $
  and   $  {\cal X}_j $ are defined in \nref{eq:dfdx}.
  In addition, we also have a multi-secant version of the no-change
  condition \eqref{eq:msecant2}.
  This is similar to a block version of the no-change condition 
  Equation~\eqref{eq:NoCh} as represented by   Equation (15) of  \cite{FangSaad07},
  which stipulates that 
  \eq{eq:NoCh2} 
  (G_j - G_{j-m}) q = 0 \quad \forall q \perp \Span \{ \FF_j \}.
  \en
 
Strong links can be established with the class of multi-secant methods to which AA belongs.
 Without loss of generality  and in an effort  to simplify notation we
also assume that $j_m = 1$ this time and $\beta=0$.
 According to
 (\ref{eq:AA}--\ref{eq:thetaj}), the $j$-th iterate becomes simply
 $x_{j+1} = x_j - \XX_j \theta_j $ where
 $\theta_j$ is a vector that minimizes $\| f_j - \FF_j \theta \|$.
 For nlTGCR(m), we have $x_{j+1} = x_j + P_j y_j $
 where $y_j $ minimizes $ \| f(x_j) + V_j y \| $. So this  is identical with
 Equation \nref{eq:AA} when $\beta = 0$  in which  $P_j \equiv \XX_j$,
 and $\FF_j\equiv V_j$:

 \bigskip 
  \begin{center}
 \begin{tabular}{c|c c c c} \hline 
   AA     &   $\XX_j$ &   $\FF_j$   & $\theta_j$         \\  \hline
   nlTGCR &    $P_j$  &  $V_j$      & $- y_j$         \\ \hline
  \end{tabular}
  \end{center}

  \bigskip 
  In  multi-secant methods we set 
  \[ \FF_j = [f_1 - f_0, f_2 - f_1, \cdots, f_j - f_{j-1}] \quad
    \XX_j = [x_1 - x_0, x_2 - x_1, \cdots, x_j - x_{j-1}] \]
  and, with  $ G_{j-m} = 0 $,    the multi-secant matrix in \nref{eq:msec} becomes
  \eq{eq:AAb0}
    G_{j} = \XX_j  (\FF_j^T\FF_j)^{-1}\FF_j^T.
  \en 
  We restate the result from \cite{FangSaad07} that characterizes multi-secant methods also known
  as Generalized Broyden
  techniques~\cite{eyert:acceleration96}, for the particular case in which
$G_{j-m} \equiv 0$:  
  \begin{itemize}
  \item
    $G_j$ in \eqref{eq:AAb0} is the only matrix
    that satisfies both the secant condition \eqref{eq:mscond} and the no-change
      condition \eqref{eq:NoCh2};
    \item $G_j$ is also the matrix that minimizes $\| G \|_F $ subject to the condition $G \FF_j = \XX_j$.
      \end{itemize} 
   
    Consider now nlTGCR. If we set $\FF_j \equiv V_j $ and $\XX_j = P_j$ in Anderson,
    the multi-secant matrix $G_j$ in \nref{eq:AAb0} simplifies to \nref{eq:msecant1M} -
    which also minimizes $\|G \|_F$ under the secant condition
    $G V_j = P_j$.
    Therefore the two methods differ mainly in the way in which the sets
    $\FF_j / V_j$ ,  and $\XX_j / P_j $ are defined. Let us use the more general notation
    $V_j, P_j$ for the pair of subspaces. 

    In both cases, a vector $v_j$ is related to the corresponding $p_j$ by the fact that
    \eq{eq:ApproxNLT} v_j \approx J(x_j) p_j . \en
    In the case of nlTGCR(m) this relation is explicitly enforced by
    a Frechet differentiation  (Line 10).     In the case of AA, we have  $v_j = \Delta f_{j-1} = f_j - f_{j-1} $ and the relation exploited
    is that
    {
    \eq{eq:ApproxAA}
    f_j \approx f_{j-1} + J(x_{j-1}) (x_j - x_{j-1}) \to
      \Delta f_{j-1} \approx J(x_{j-1}) \Delta  x_{j-1} .
      \en } 

      However, relation \nref{eq:ApproxNLT} in nlTGCR is \emph{more accurate}
      because we use an additional function evaluation to explicitly obtain an
      accurate approximation (ideally exact value) for $J(x_j) r_j$ in line
      10 {which is likely to lead to a smaller error in \nref{eq:ApproxNLT}.} In
      contrast when $x_j$ and $x_{j-1}$ are not close, then \nref{eq:ApproxAA}
      can be a very rough approximation. This is a key difference between the
      two methods.

\subsection{Line search techniques and convergence analysis}
\label{subsec:linesearch}
In this section, we discuss how to include  line search
techniques to improve the global convergence of nlTGCR(m). More specifically, if
$d_j=P_jy_j$, then Line 7 of Algorithm \ref{alg:nltgcr} will be replaced by
$x_{j+1} = x_j + \alpha_j d_j$ with a suitable stepsize $\alpha_j$.  

When nlTGCR(m) is applied to solve a nonlinear system $f(x) = 0$, we define the
scalar function $\phi(x) = \frac{1}{2}\|f(x)\|_2^2$.  Since
$\nabla\phi(x) = J(x)^T f(x)$ and $r_j = -f(x_j)$, the stepsize $\alpha_j$ is
chosen to fulfill the Armijo-Goldstein condition~\cite{armijo1966minimization}:
\begin{equation}
	\|f(x_j + \alpha_j d_j)\|_2^2 \leq \|r_j\|_2^2 - 2c_1\cdot \alpha_j \langle J(x_j)^{T}r_j, d_j \rangle.
	\label{eq:Armijo3}
\end{equation}
If we approximate the nonlinear residual $-f(x_{j+1})$ with the linearized one
$ r_{j+1} = r_j - V_j y_j$, we need to replace the left hand side of the inequality
\eqref{eq:Armijo3} by $\|r_j - \alpha_j V_j y_j\|^2_2$. Here $d_j$ is a descent
direction for $\phi$ if $\langle J(x_j)^Tr_{j},d_j \rangle > 0$.

It is possible to inexpensively  check the above condition by first noting that
\[
\langle J(x_j)^T r_j,d_j \rangle
=
\langle  r_j, J(x_j) d_j \rangle .
\]
Then,  using Frechet differentiation we can write 
\eq{eq:anglChk}
\langle J(x_j)^T r_j,d_j \rangle
\approx 
\frac{1}{\eps} 
\langle  r_j, f(x_j + \eps d_j) - f(x_j) \rangle
\en
where $\eps$ is some small parameter.
We already have $f(x_j)$ available and in the context of line search techniques
$f(x_j+\lambda_0 d_j)$ is computed as the first step, where $\lambda_0$ is some
parameter, often set to $\lambda_0 =1$.
If $\| \lambda_0 d_j \|$ is small enough relative to $x_j$  we can get
an accurate estimate of
$\langle J(x_j)^T r_j,d_j \rangle$ using relation \nref{eq:anglChk} with
$\eps$ replaced by $\lambda_0$.
If not, we may compute an additional function evaluation to obtain
$f(x_j + \eps d_j) $ for a small $\eps $ to get the same result with high accuracy.
Note that as the process nears convergence
$\| d_j\|_2 $ becomes small and this is unlikely to be needed.

Often, the direction $d_j = P_jy_j$ is observed to be a descent direction and
this can be explained from the result of Proposition~\ref{prop:nltgcr}.
A
few algebraic manipulations lead to the following proposition.

\begin{proposition}\label{prop:3}
	Let $\phi(x) = \half \| f(x)\|_2^2 $ and let
	$\tilde v_{j_m}, \cdots, \tilde v_j$ be the columns of:
	\eq{eq:tilV}
	\tilde V_j \equiv J(x_j) P_j . \
	\en 
	Then, 
	\eq{eq:prop3} 
	\langle \nabla \phi(x_j), d_j\rangle
	{=  -  \sum_{i=j_m}^{j} \langle v_i, r_j \rangle  \langle \tilde v_{i}, r_{j} \rangle
        = - r_j^T \tilde V_j V_j^T r_j .} 
	\en
\end{proposition}
{
\begin{proof} 
	From  Property~(4) of Proposition~\ref{prop:nltgcr} we have
	$  V_j^T r_j = \langle v_{j}, \tilde r_{j}\rangle e_{m_j} -  V_{j}^T z_{j} $ and so
\begin{align} 
\langle\nabla \phi(x_j), d_j\rangle
  &= - \langle J(x_j)^T r_j, P_j [ \langle v_{j}, \tilde r_{j}\rangle e_{m_j} -  V_{j}^T z_{j} ]\rangle
  \nonumber \\
  &= - \langle v_{j}, \tilde r_{j}\rangle \langle J(x_j)^T r_j, p_j\rangle + \langle J(x_j)^T r_j, P_jV_{j}^T z_{j}\rangle
  \nonumber \\
  &= - \langle v_{j}, \tilde r_{j}\rangle \langle r_j, J(x_j) p_j\rangle + \langle r_j,  J(x_j) P_j V_{j}^T z_{j} \rangle
  \nonumber \\ 
  &= - \langle v_{j}, r_j+z_j\rangle \langle \tilde v_{j}, r_{j}\rangle  + \langle r_j, \tilde V_j V_{j}^T z_{j} \rangle
  \nonumber \\
  &= - \langle v_{j}, r_j\rangle \langle \tilde v_{j}, r_j\rangle
    -  \langle v_{j}, z_j\rangle \langle\tilde v_{j}, r_{j}\rangle + \langle r_j, \tilde V_j V_{j}^T z_{j} \rangle .
		\label{eq:PrfDesc1} 
\end{align}
We write the last term 
$ \langle r_j, \tilde V_j V_j^T z_j \rangle$ in the form
\[ \sum_{i=j_m}^j  r_j^T \tilde v_i v_i^T z_j = \sum_{i=j_m}^j  \langle\tilde v_i, r_j\rangle \langle v_i, z_j\rangle . \]
 We now note that the last term in the above sum is equal to the middle term in \nref{eq:PrfDesc1}.
 The result is that 
 \eq{eq:PrfDesc2} 
 \langle \nabla \phi(x_j), d_j\rangle
 = - \langle v_{j}, r_j\rangle  \langle \tilde v_{j}, r_j\rangle
 +  \sum_{i=j_m}^{j-1} \langle r_j, \tilde v_i v_{i}^T z_{j}\rangle.
 \en
 Next observe that for each of the terms in the sum we have:
 \[
   v_{i}^T z_{j} = v_{i}^T (\tilde r_j - r_j) = -v_{i}^T r_j \to
    \sum_{i=j_m}^{j-1} \langle r_j, \tilde v_i v_{i}^T z_{j}\rangle 
= - \sum_{i=j_m}^{j-1}  \langle\tilde v_i, r_j\rangle \langle v_i, r_j\rangle .
\]
Substituting this in \eqref{eq:PrfDesc2} yields the desired result.
\end{proof}
}

{We have  $v_i \approx J(x_i) p_i$ while $\tilde v_i  = J(x_j) p_i$.}
Near convergence, the two vectors will be close enough that the
inner products
$\langle v_i, r_j \rangle $ and $ \langle \tilde v_{i}, r_j \rangle $ will have the same sign in which case
the inner product \nref{eq:prop3} is negative and $d_j$ is a descent
direction.
In addition, we saw in the proof that $\langle v_i, r_j \rangle = -\langle v_i, z_j \rangle $ for $j_m \le i \le j-1$.
{Thus, the terms $ \langle v_i, r_j\rangle \langle \tilde v_{i}, r_{j}\rangle $ are likely to be smaller order terms for
$i<j$.}

{Based on the above discussion, we may assume that the direction produced by
Algorithm \ref{alg:nltgcr} is likely to always be a descent direction but this not guaranteed.
However, if needed, this can
be inexpensively checked as described earlier at a small additional cost.}

 In order to ensure global convergence, we also implemented backtracking line
 search. Specifically, an initial stepsize $\alpha_j^{(0)}$ is defined where the
 superscript indicates the backtracking steps. We repeatedly set
 $\alpha_j^{(k+1)} := \tau\cdot\alpha_j^{(k)}$ for a shrinking parameter
 $\tau\in(0,1)$ and check if the Armijo-Goldstein condition is satisfied or if
 we have exceeded the maximum allowed number of backtracking steps. Our tests
 use $\tau = 0.8$. In addition, we found it effective to select the initial
 stepsize $\alpha_j^{(0)}$ adaptively in order to reduce the number of line
 search steps. For example, letting $\alpha_0^{(0)} := 1$, we define
\begin{equation}
	\alpha_{j+1}^{(0)} := 
	\left\{
	\begin{array}{ll}
		\min\{1,~\alpha_j^{(0)}/\tau\}, & \text{if line search finishes in 1 step;} \\
		\tau\cdot\alpha_j^{(0)}, & \text{otherwise.}
	\end{array}\right.
\end{equation}

It is possible to prove that under a few assumptions nlTGCR(m) converges
globally.  Suppose that we have a line search procedure
which at the $j$th step considers iterates of the form \eq{eq:ls1}
x_{j+1} (t) = x_j + t d_j, \quad \mbox{with} \quad d_j = P_j y_j = P_j
V_j^T r_j .  \en We further assume that the line search is exact, i.e.,
  
  \bigskip\noindent\textbf{Assumption~1:}
  \[
    x_{j+1} = \text{argmin}_{x_{j+1}(t), t  \in [0, 1]} \| f(x_{j+1}(t)) \| . \]
  We expand $f(x_j + t d_j)$ locally as follows
 \eq{eq:ls0}
 f(x_j+t d_j) = f(x_j) + t  J(x_j) d_j +  s_j(t) . 
 \en

 The term $s_j(t) $ is a second order term and  we will make the following
 smoothness  assumption:

 \bigskip\noindent\textbf{Assumption~2:}  There is a constant $K>0$ such that for each $j$ we have:
 \eq{eq:asp1} 
   \| s_j(t) \| \le  K \| f(x_j) \|  t^2
   \en
In addition, we will assume, without loss of generality, that $K \ge \frac12$. 

   Next, for the term $J(x_j) d_j$ in \nref{eq:ls0},
   which is equal to $J(x_j) P_j y_j$, it is helpful to
   exploit the notation \nref{eq:wi}  and observation \nref{eq:Jxjpi}.
   Indeed, these imply that:
   \eq{eq:eqwj} J(x_j) P_j y_j = V_j y_j + W_j y_j  . \en
   It was argued in the discussion following Proposition~\ref{prop:2}  that $W_j y_j$
   can be expected to be much 
   smaller in magnitude than $r_j$ and this leads us to our 3rd assumption.

 \bigskip\noindent\textbf{Assumption~3:} There is a scalar $0\le \mu<1$ such that for all $j$: 
 \eq{eq:eta}     \| W_j y_j \| \le \mu \| r_j\| \ .       \en 
   
 Finally, we assume that each linear least-squares problems is solved with a certain relative tolerance.

 \bigskip\noindent\textbf{Assumption~4:}  At every step $j$ the least-squares solution $V_j y_j$ satisifies
   \eq{eq:lemAp2} 
     \| f(x_j) + V_j y_j \| \le \eta \ \| f(x_j) \| . 
     \en

     With these the following theorem can be stated.
     \begin{theorem}\label{thm:conv}
       Let Assumptions 1--4 be satisfied
     and assume that  the linear least-squares problem in Eq. \nref{eq:lemAp2}
       is solved with a relative tolerance $\eta$
       with $0<\eta < 1-\mu$, i.e., $\eta $  satisfies
       \eq{eq:ThAs1}  c \equiv 1 - (\eta+\mu) >0 . \en 
       Then, under these assumptions, we have $0 < 1 - c^2/4K < 1 $ and
       the following inequality is satisfied at each step $j$:
       \eq{eq:main}
    \| f(x_{j+1}) \| \le \left[1 - \frac{c^2}{4K} \right]  \ \| f(x_{j}) \| .
       \en
       \end{theorem} 
       
       \begin{proof}
Recall that
 \eq{eq:ls2}
 f(x_{j+1}(t))  = f(x_j) + t  V_j y_j + t W_j y_j + s_j(t)
 \en
which we rewrite as: 
   \eq{eq:ls3}
 f(x_{j+1}(t)) = t [f(x_j) + V_j y_j] + (1-t)f(x_j)  +   t W_j y_j + s_j(t) . 
 \en
  Thus, for any $t$ with  $0\le t\le 1$ 
  \begin{align} 
   \|f(x_{j+1}(t)) \| & \le  t \| f(x_j) + V_j y_j\| + [(1-t)+\mu t] \| f(x_j) \|
                        + K \| f(x_j) \| t^2 \nonumber \\
   & \le  t \eta \| f(x_j) \|  + [(1-t)+\mu t] \| f(x_j) \|  + K \| f(x_j) \| t^2  \nonumber \\
   & \le \left[ 1 - (1 - \eta -\mu) t \right] \| f(x_j) \|  + K \| f(x_j)\|  t^2  .  \label{eq:ls5} 
     \end{align}
     Let us set $\rho_j = \| r_j\|$.
     From the definition of $c$ and  \nref{eq:ls5}  we get
     \eq{eq:ls7}
 \min_t \|f(x_{j+1}(t)) \| \le  \rho_j \times \min_t \left[  K t^2 - c t +1 \right] . 
 \en 
 The minimum with respect to $t$ of the quadratic function in the brackets  is reached for
 $t_{opt}  = c /(2K)$ and the minimum value is 
 \eq{eq:ls6} 
   \min_t 
    \left[K t^2 - c t + 1 \right]  =
   \left[ 1 - \frac{c^2}{4 K} \right]  . 
 \en
 Note that our assumptions imply that $c > 0$ so $t_{opt} >0$.
 In addition, $c\le 1$ and so 
 $c/(2 K) \le 1/(2K) \le 1$ (since $K \ge 1/2$). Thus,
 the (exact) line search for $t \in \ [0, 1] $ will yield $x_{j+1}$. 
 Relations~\nref{eq:ls6} and \nref{eq:ls7} imply that $\rho_{j+1} \le [1-c^2/(4K)] \rho_j$
 where $ | 1-c^2/(4K)| < 1 $ under the assumptions on $c$ and $K$. This proves
 relation \nref{eq:main} which establishes convergence under the stated assumptions. 
\end{proof}

The theorem suggests 
that the residual norm for each solve must be reduced by a certain minimum amount in order
 to achieve convergence.
 For example if $\mu = 0.1$ (recall that $\mu $ is small under the right
 conditions) and we set $c = 0.1$   then we would need  a reduction of at least 0.8, {which is
 similar to the residual norm reduction required in our experiments by default. Also, the four assumptions do not require $f$ to be convex.}

\subsection{Extension to the stochastic case}
{
Inspired by the great success of nonlinear acceleration methods in accelerating fixed-point iterations, it is natural to ask whether they can be applied to accelerate practical applications, such as stochastic optimization problems \cite{Robbins1951ASA} in deep learning where gradients are approximated by using a random batch of samples in each iteration.
We show that the nlTGCR(m)} algorithm can be easily generalized to stochastic cases such as
mini-batching in deep learning. The main difference is that we now build the
pair $\{P_j, V_j\}$ with respect to different objective functions $\phi_j(x)$ instead of {only one objective function $\phi(x)$. In this case, Line 10 of Algorithm \ref{alg:nltgcr} becomes $v := J_j(x_{j+1}) p$
where $J_j(x)$ indicates the Jacobian corresponding to $\phi_j(x)$.} We found it important to also add gradient pre-normalization to enhance the
stability and speed up the training process.
That is, the combined gradient of all layers of the entire model has unit 2-norm.
{Pre-normalizing the gradient is
useful because the ``search'' space $\operatorname{Span}\{P_j\}$ is invariant to the
magnitude of its component vectors $p_j$.} 
When the batch size is sufficiently
large, the ``search'' direction originating from the stochastic gradient is
close to the one generated from the full-batch gradient.
Pre-normalization can help mitigate the damaging impact of  
noisy gradients\cite{cabana2021backward} 
{Given the significance of this topic, we intend to explore this in future work, focusing on the theoretical aspects of the proposed methods}
\footnote{Sharath Sreenivas,
	Swetha Mandava, Boris Ginsburg and Chris Forster. Pretraining BERT with
	Layer-wise Adaptive Learning Rates. NVIDIA, December 2019.}.

\section{Numerical experiments}
This section presents a few experiments to compare  nlTGCR(m) with existing methods in the
literature.  We propose an adaptive mechanism that combines the
nonlinear update version with the linearized residual computation.  This
adaptive version is presented in Section \ref{subsec:ada} and is implemented in other
experiments by default. All experiments were conducted on a workstation
equipped with an Intel i7-9700k CPU and an NVIDIA GeForce RTX 3090 GPU. The
first three experiments were implemented in {MATLAB} 2022b, and the baseline
methods were based on the implementations by H. Kasai~\cite{gdlibrary} and
C.T. Kelley~\cite{kelley_solving_2003}. Deep learning experiments are implemented using PyTorch \cite{pytorch} and run with GPU acceleration. 

\subsection{Bratu problem}
\label{sec: bratu}
We first consider a nonlinear Partial Differential Equation (PDE) problem, namely the Bratu
problem~\cite{hajipour2018accurate}  of the following form with  $\lambda = 0.5$: 
\begin{align*} 
   \Delta u + \lambda e^ u &= 0  \ \ \text{in} \ \ \Omega = (0, 1)\times (0,1)  \\
   u(x,y) &= 0 \ \text{for} \
            (x,y) \in \partial
            \Omega  
  \end{align*} 
 This problem is known to be not particularly difficult to solve but our purpose
 is to illustrate the importance of an accelerator that exploits symmetry.
 The problem is discretized with  Centered Finite
 Differences~\cite{chaudhry2008open,folland1995introduction,wilmott1995mathematics}
 using a grid of size $100 \times 100 $, yielding a 
 nonlinear system of equations $f(s) = 0$ where 
 $f:\mathbb{R}^{n}\to\mathbb{R}^{n}$, with $n=10,000$. The corresponding fixed point
 problem takes the form $g(s) = s + \beta f(s)$.
	
\subsubsection{Superiority of the adaptive update version}
\label{subsec:ada}

The Bratu problem possesses the property of being almost linear despite the presence
of an exponential term. This property also appears in numerous applications when
nearing convergence. 
{
Hence, this problem is suitable to illustrate the difference between the linearized and nonlinear update versions, while also showing the advantage of utilizing the adaptive update version.
}

\begin{figure}[htb]
 	\centering
 	\subfigure[Number of function evaluations vs. relative residual norm.]{
 		\includegraphics[width=0.46\textwidth]{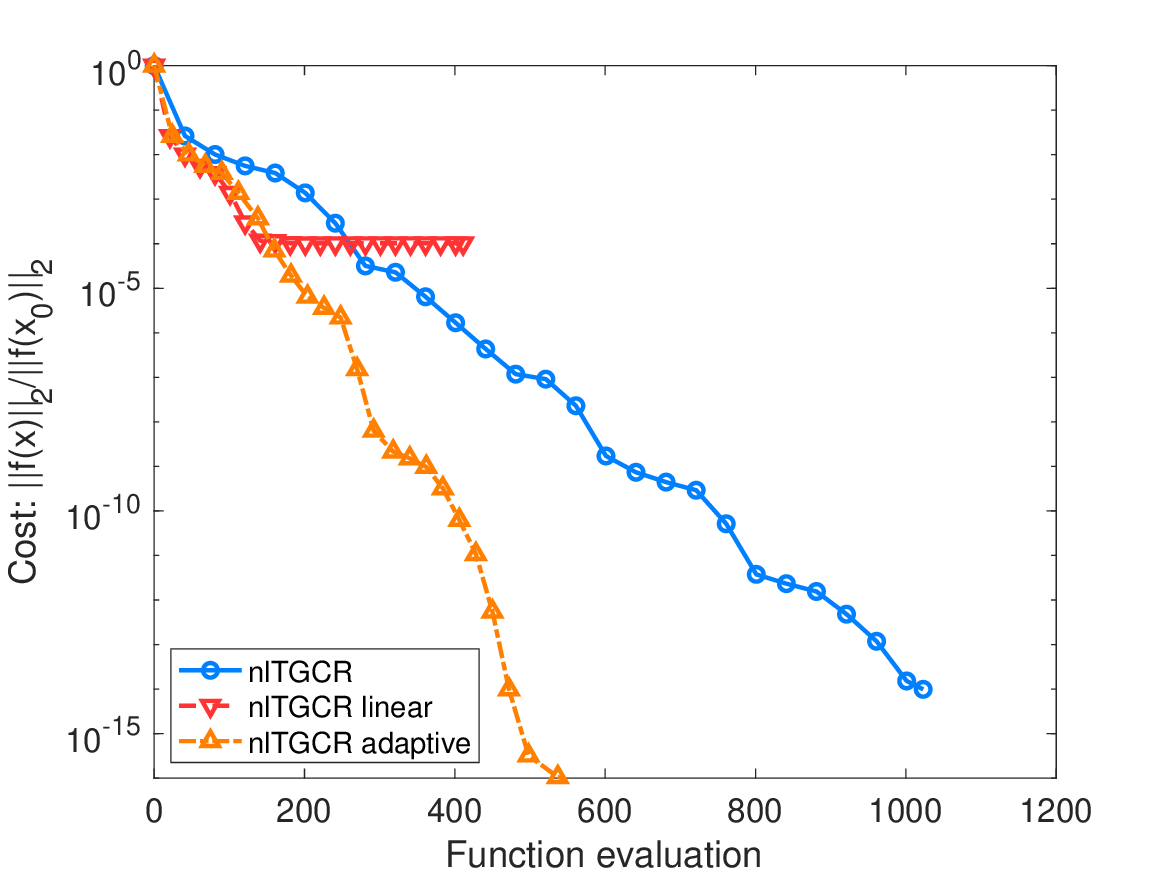}
 		\label{fig:bratu_mode_1}
 	}
 	\subfigure[Number of iterations vs. relative residual norm.]{
 		\includegraphics[width=0.46\textwidth]{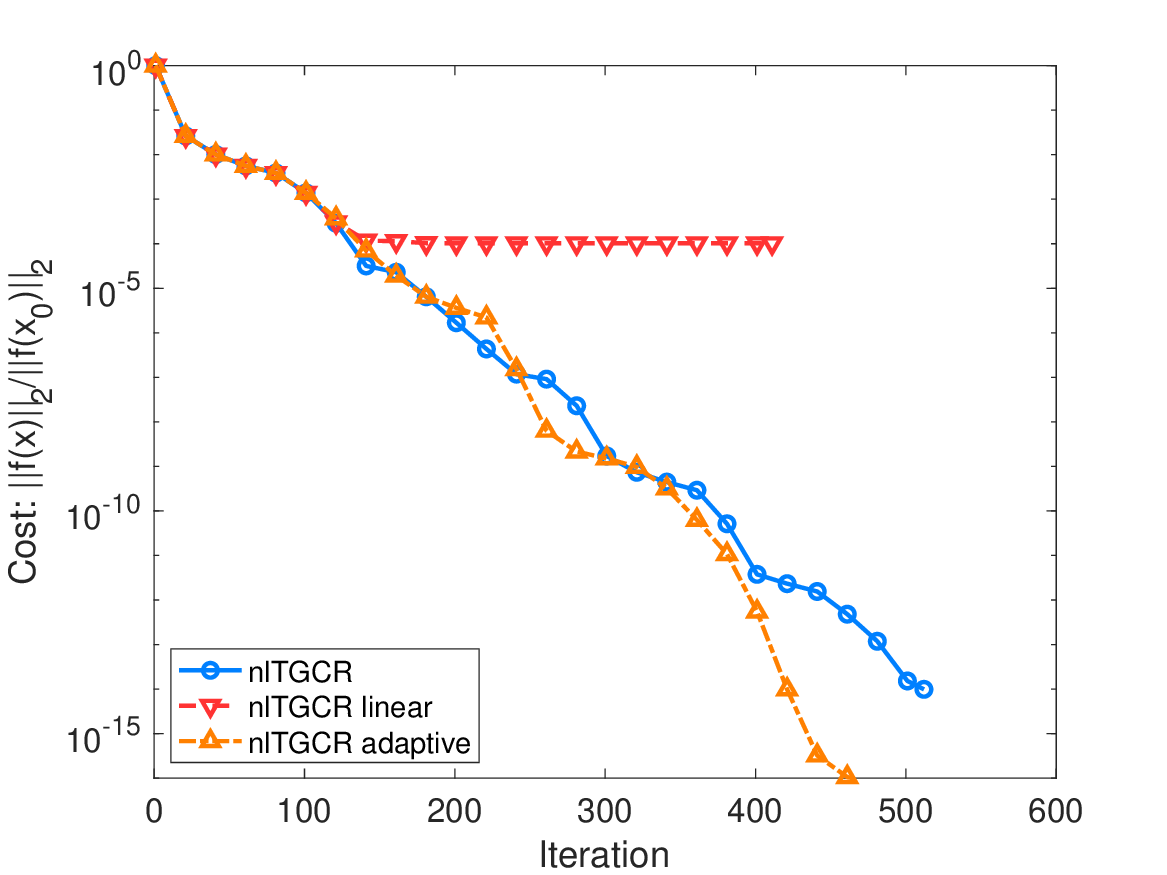}
 		\label{fig:bratu_mode_2}
 	}
 	\caption{Comparison between the standard, linearized update, and adaptive update versions of nlTGCR(m) with $m=1$ on the Bratu problem. Each marker represents 20 iterations.}
  \label{fig:bratu_mode}
\end{figure}

For the experiment in Figure \ref{fig:bratu_mode}, the window size is $m=1$, and
the starting point is a vector of all-ones $x_0 = \textbf{1}$. The
cost/objective is the relative residual norm $\|f(x)\|_2/\|f(x_0)\|_2$.  We
compare all three types of residual update schemes in terms of the number of
function evaluations and present results in Figure \ref{fig:bratu_mode_1}. It
can be observed that the convergence rate of the adaptive update version is
close to that of the linearized updated version in the first few
iterations. This is because the adaptive update version switches to linear
updates at the second iteration and switches back to the nonlinear form at the
110th iteration where the linear update version stalls. As for the cost (Y-axis)
of each iteration, Figure \ref{fig:bratu_mode_2} indicates that all three
versions decrease almost in the same way per iteration before the onset of
stagnation. The adaptive update version performs as well as the nonlinear update
version afterward. {Sections \ref{subsec:sym} and \ref{sec:molecular} report on more
experiments with the adaptive version of nlTGCR(m). In Sections
\ref{sec:resnet} to \ref{sec:gcn} we utilize the standard (nonlinear) update
version of nlTGCR(m) for deep learning tasks, as the proposed residual check is
not applicable in a stochastic context.}
	
\subsubsection{Exploiting symmetry}
\label{subsec:sym}
\begin{figure}[htb]
 	\centering
 	\subfigure[Starting point $x_0 = \textbf{0}$.]{
 		\includegraphics[width=0.46\textwidth]{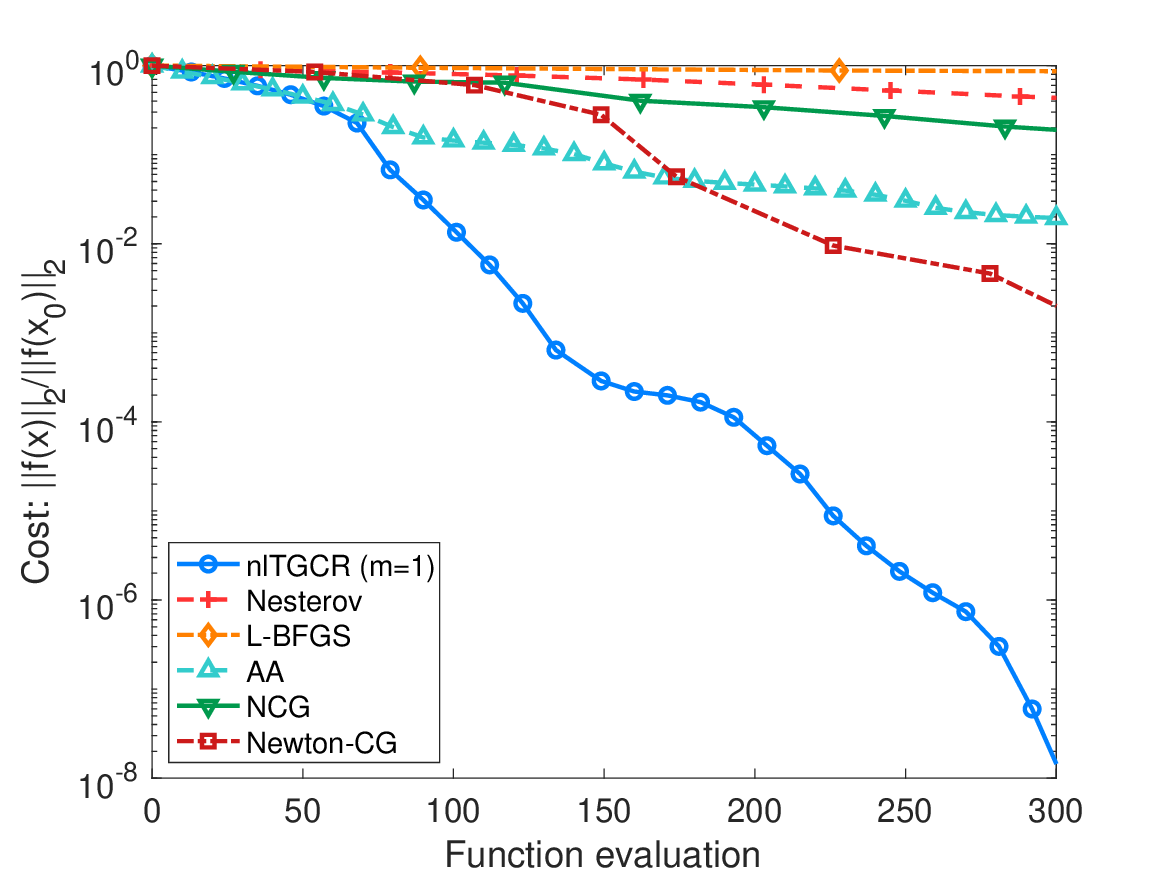}
 		\label{fig:bratu_1}
 	}
 	\subfigure[Starting point $x_0 = \textbf{1}$.]{
 		\includegraphics[width=0.46\textwidth]{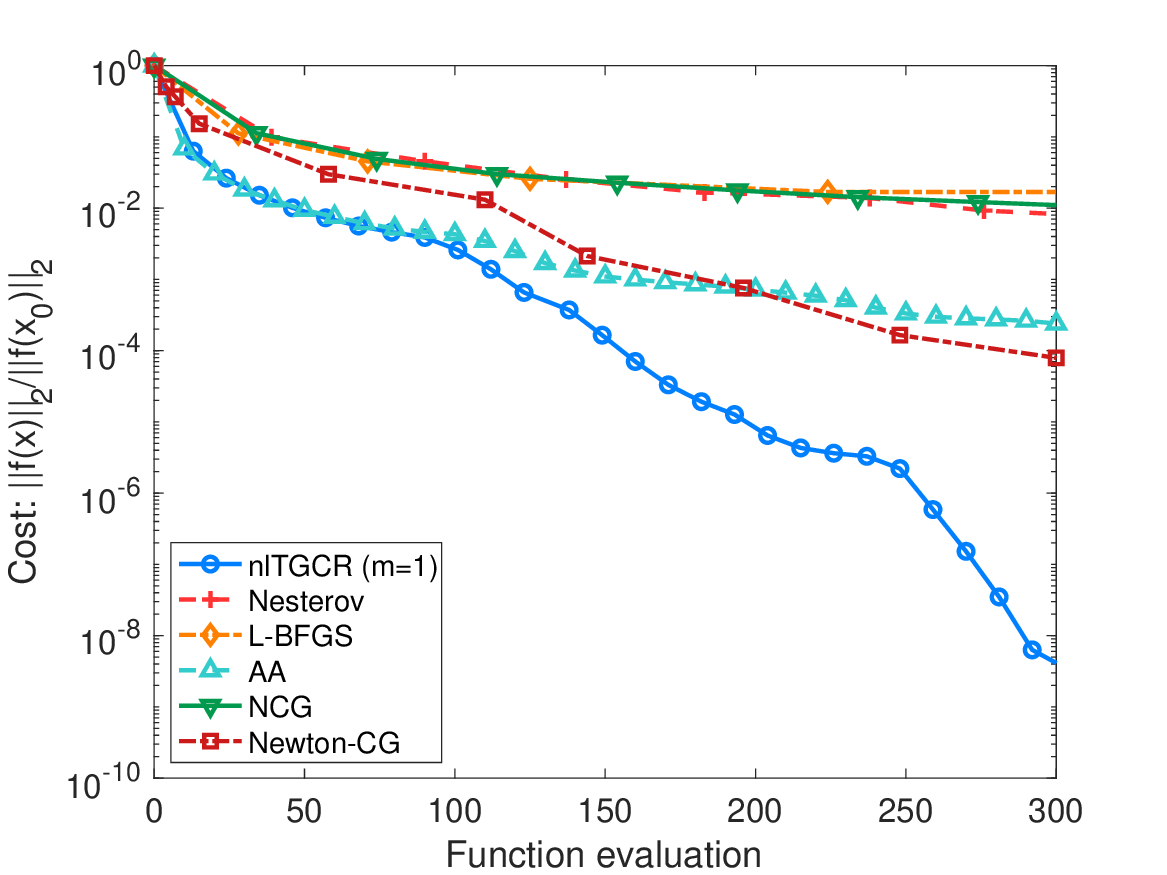}
 		\label{fig:bratu_2}
 	}
 	\caption{Number of function evaluations vs. relative residual norm on the Bratu problem with different starting points. Each marker represents 10 iterations except Newton-CG where each marker represents 1 outer loop step.}
\end{figure}

We now investigate whether nlTGCR(m) takes advantage of short-term recurrences when
the Jacobian is symmetric.
Recall that in the linear case,  GCR is mathematically equivalent to
the Conjugate Residual (CR) algorithm when the matrix is  symmetric. In this
situation a  window size $m=1$ is optimal~\cite{Saad-Book2}.  We tested
nlTGCR(m) along with baselines including Nesterov's Accelerated
Gradient~\cite{nesterov}, L-BFGS~\cite{liu1989limited}, AA, Nonlinear
Conjugate Gradient (NCG) of Fletcher--Reeves'
type~\cite{fletcher1964function}, and Newton-CG~\cite{Dembo-al}. Results
are presented in Figures \ref{fig:bratu_1} and \ref{fig:bratu_2}. We
analyzed the convergence in terms of the number of function evaluations
rather than the number of iterations because the backtracking line
search is implemented for all methods considered except  AA by
default. We present the results of the first 300 function evaluations for all methods. The window size for nlTGCR(m) is $m=1$, while for L-BFGS and AA, it is $m=10$. The
mixing parameter for AA is set to $\beta = 0.1$ as suggested in
\cite{brezinski2022shanks}. For Newton-CG method, the maximum number of steps in the
inner CG solve is 50. This inner loop can be terminated early if
\begin{equation}
  \|r_k\| \leq \eta \|f(x_0)\|.
\end{equation}
We choose the forcing term $\eta = 0.9$ and adjust it via
the Eisenstat-Walker method~\cite{Eisenstat-Walker}.

	The Jacobian for the Bratu problem $\nabla f(s)$ is symmetric positive
        definite (SPD), making nlTGCR(m) with $m=1$ a highly efficient
        method. Other methods that do not take advantage of this symmetry
        require a larger number of vectors to achieve comparable performance. We
        tested the methods with two different initial guesses. The first,
        used in the experiments in Figure \ref{fig:bratu_1}, is
        $x_0 = \textbf{0}$ which  is somewhat close to the global optimum. The
        second initial guess, used in the experiments in Figure
        \ref{fig:bratu_2}, is the  vector of all ones $x_0 = \textbf{1}$. In both
        cases, we set the window size of L-BFGS and AA to 10, which means  that 20
        vectors in addition to $x_j$ and $r_j$ need to be stored. In contrast, nlTGCR(m=1)
        only requires 2 extra vectors. In this problem, nlTGCR(m=1), Nesterov, NCG, and Newton-CG are memory-efficient in terms of the number of vectors required to store information from previous steps. Among these competitive
        methods, nlTGCR($m=1$) performs best -- suggesting that
        nlTGCR(m) benefits from  symmetry.

{
        To further clarify our claim, we modify the Bratu problem:
\begin{equation*} 
   \Delta u + \alpha u_x + \lambda e^ u = 0
\end{equation*} 
by introducing asymmetry to the Jacobian via the term $\alpha u_x$. We employ the standard nlTGCR method with no restarts, and we consider various table sizes $m = 1, 2, 3, 5$, 10 and initiate the process using an all-one vector $\mathbf{1}$.

\begin{figure}[htb]
 	\centering
 	\subfigure[Near-linear symmetric Jacobian.]{
 		\includegraphics[width=0.46\textwidth]{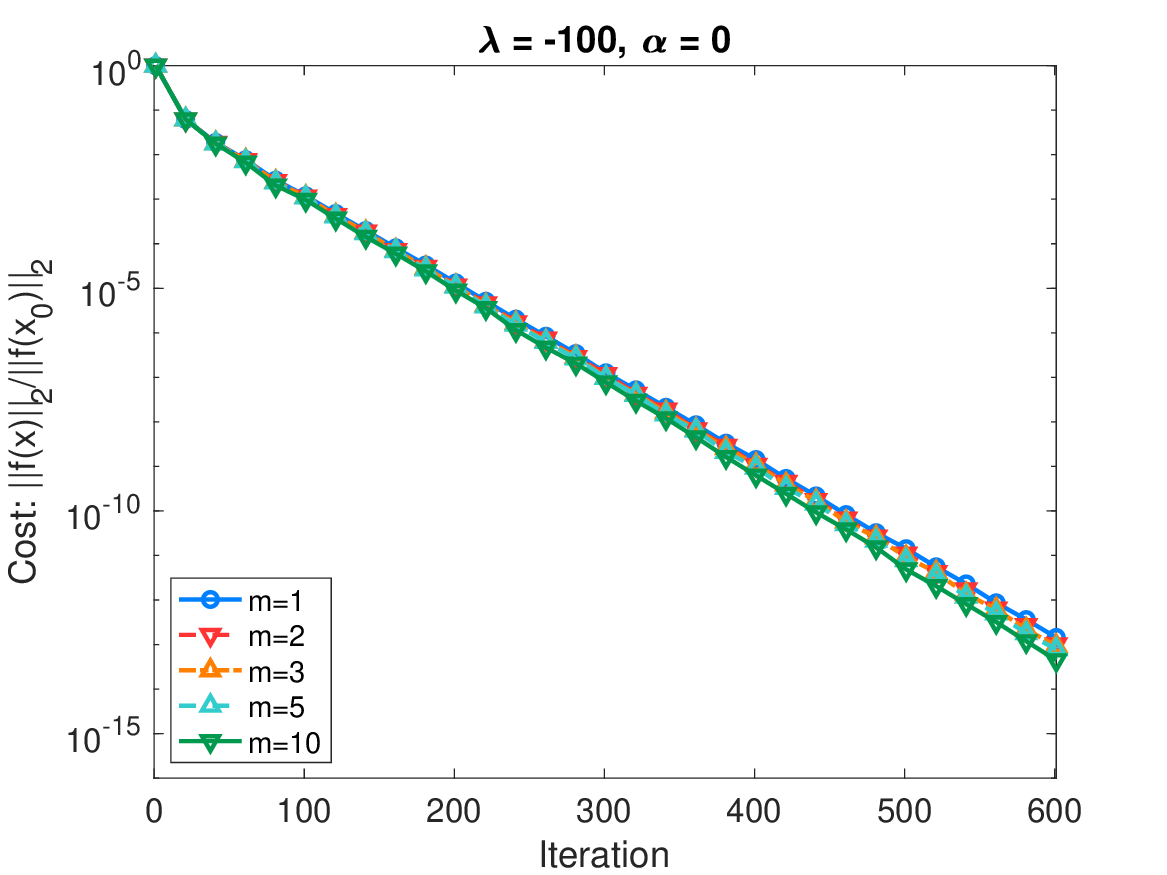}
 		\label{fig:mod_bratu_1}
 	}
 	\subfigure[Linear asymmetric Jacobian.]{
 		\includegraphics[width=0.46\textwidth]{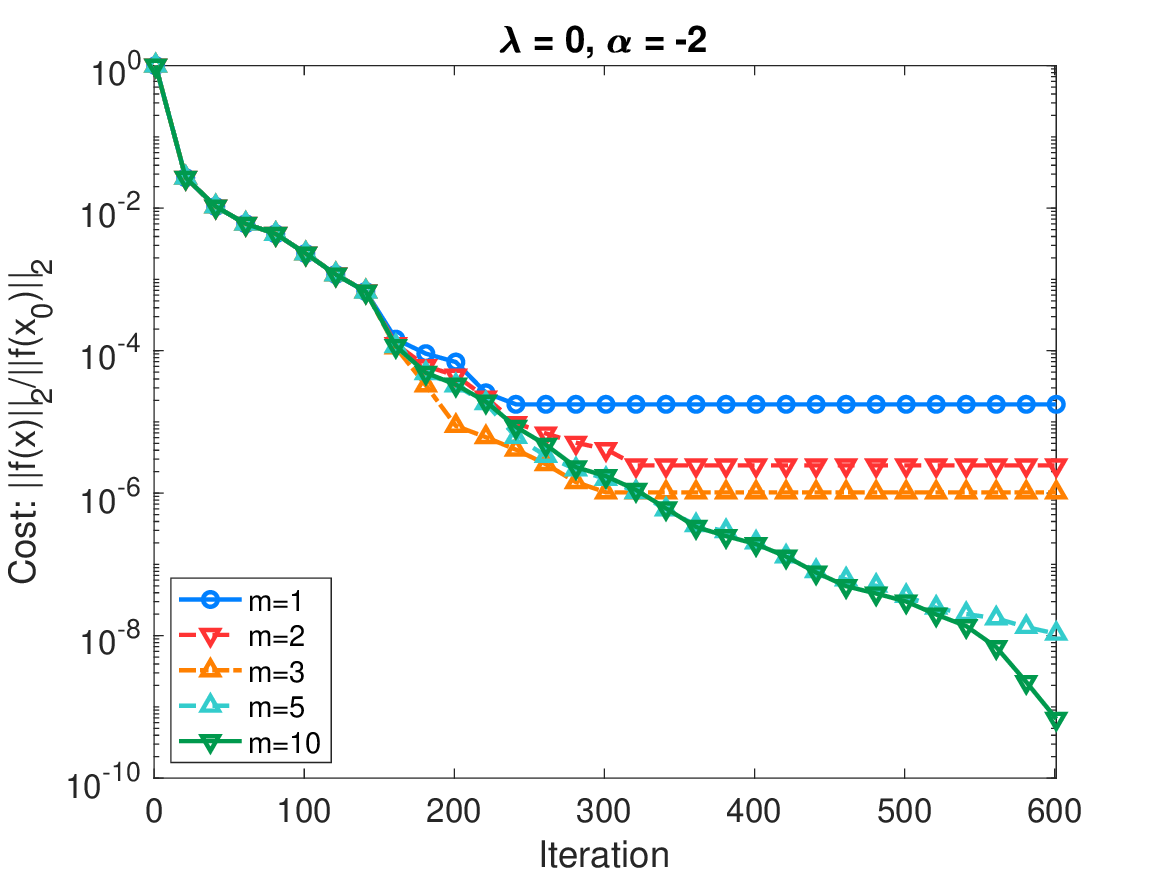}
 		\label{fig:mod_bratu_2}
 	}
 	\caption{Comparison of nlTGCR(m) with $m=1,2,3,5,10$ on the modified Bratu problem. Each marker represents 20 iterations.}
  \label{fig:bratu}
\end{figure}

In Figure \ref{fig:mod_bratu_1}, we set $\lambda=-100$ and $\alpha=0$, resulting in a nonlinear symmetric Jacobian. The overlap of all lines confirms the observation that nlTGCR benefits from the symmetry of the Jacobian. Conversely, in Figure \ref{fig:mod_bratu_2} where $\lambda = 0$ and $\alpha = -2$, the Jacobian is linear and asymmetric. The selection of table size $m$ profoundly impacts the performance. Hence, nlTGCR exploits the symmetry to establish short-term recurrence and improve convergence, while the nonlinearity of the problem mainly affects the adaptive update mechanism.
}
	
	\subsection{Molecular optimization with Lennard-Jones potential}
	\label{sec:molecular}
\begin{figure}[htbp]
   \centering
\subfigure[Geometrical configurations.]{
     \includegraphics[width=0.46\textwidth]{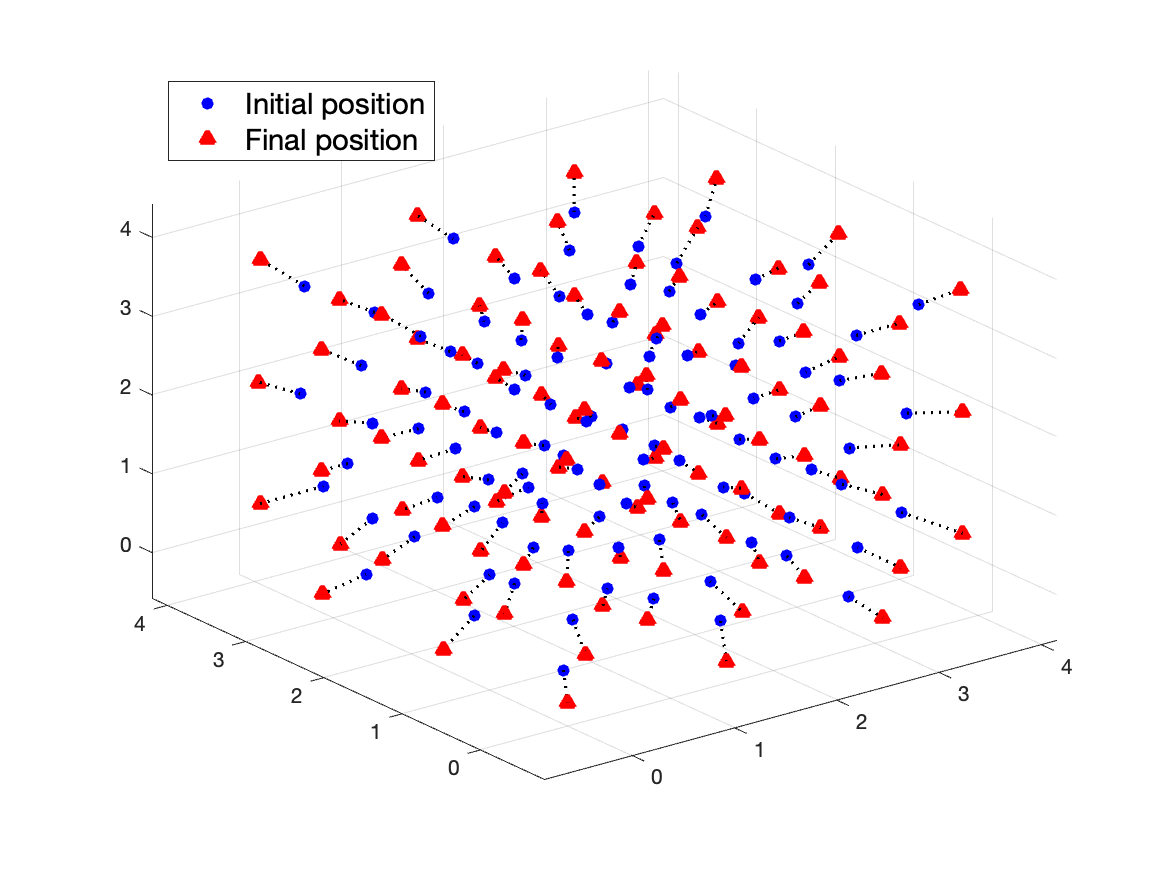}
     \label{fig:LJ_2}   }
   \subfigure[Convergence of various methods.]{
     \includegraphics[width=0.46\textwidth]{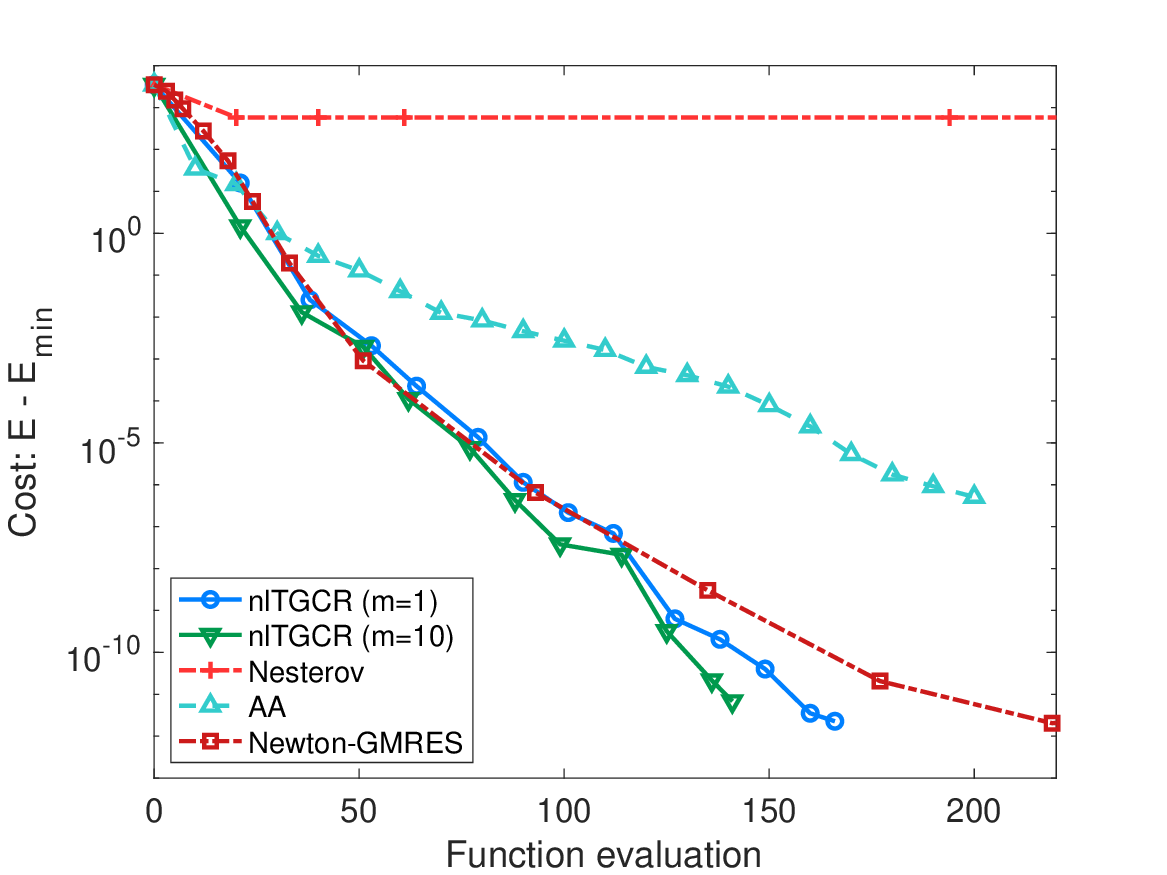}
     \label{fig:LJ_1}
   }
   \caption{(a) Initial and final configurations of 108 atoms with the Argon cluster experiment. (b) Number of function evaluations
vs. shifted potential norm on the Lennard-Jones problem. Each marker represents 10 iterations for all methods except Newton-GMRES where each marker represents 1 outer loop step.}
   \label{fig:LJ}
\end{figure}

	The second experiment focuses on the molecular optimization with the Lennard-Jones (LJ) potential which is a geometry optimization problem. The goal is to find atom positions that minimize total potential energy as described by the LJ potential\footnote{Thanks: We benefited from Stefan Goedecker's course site at  Basel University.}:
	\begin{equation}
		E = \sum_{i=1}^{N} \sum_{j=1}^{i-1} 4 \times
		\left[\frac{1}{\|y_i - y_j\|^{12}} - \frac{1}{\|y_i - y_j\|^{6}} \right].
	\end{equation}
	In the above expression each $y_i$  is a 3-dimensional vector whose components are the coordinates of the location of
	atom $i$. We start with a certain configuration and 
	then optimize the geometry by minimizing the potential starting from that position. Note that the resulting position
	is not a global optimum - it is just a local minimum around the initial configuration (see, e.g., Figure \ref{fig:LJ_2}). In this particular example,
	we simulate an Argon cluster by taking
	the initial position of the atoms to consist of a perturbed initial Face-Centered-Cubic (FCC) structure \cite{meyer1964new}.
	We took 3 cells per direction - resulting in 27 unit cells. FCC cells include 4 atoms each and so we end up with a
	total of 108 atoms. 
	The problem is rather hard due to  the high powers in the potential. In this situation  backtracking or some form of
	line search is essential.

        In this experiment, we set $f=\nabla E$.  Each iterate in nlTGCR(m) is a vectorized array of
        coordinates of all atoms put together. So, it is a flat vector of length
        $3\times 108 = 324$.  We present the results of the first 220 function evaluations for nlTGCR(m), AA, Nesterov, and Newton-GMRES in Figure \ref{fig:LJ_1}. The
        reason for excluding L-BFGS, NCG, and Newton-CG is that the Jacobian/Hessian
        of the LJ problem is indefinite at some $x_j$ which can lead to a
        non-descending update direction. The window sizes for nlTGCR(m) are $m=1$ and $m=10$, while for AA it is $m=10$ and for Newton-GMRES it is $m=20$. This is because nlTGCR(m) and AA require storing twice as many additional vectors as the window size to generate the searching subspace, while Newton-GMRES does not. For each inner GMRES solve, GMRES is allowed to run up to 40 steps and utilizes a forcing term $\eta=0.9$ to terminate the inner loop. Moreover, AA does not converge unless the underlying fixed-point iteration $s_{j+1} = s_j + \beta f(s_j)$ is carefully chosen. In this experiment, we select $\beta = 10^{-3}$. The cost (Y-axis) represents the shifted potential $E - E_{min}$, where $E_{min}$ is the minimal potential achieved by all methods, approximately equal to $-579.4638$.

        In Figure \ref{fig:LJ_1}, we observe that nlTGCR($m=10$) converges the fastest.
        The convergence of nlTGCR($m=1$) and Newton-GMRES with a subspace dimension of 20 is close and
slightly slower than nlTGCR($m=10$). One observed phenomenon worth mentioning is the quick
termination of the inner loop of Newton-GMRES during the first few
iterations. Newton-GMRES moves quickly to the next Jacobian  at the beginning and
focuses on a single Jacobian when nearing convergence. This behavior is
made possible by the use of the Eisenstat-Walker technique, which accounts
for the fast convergence of Newton-GMRES. However, without this early stopping mechanism, Newton-GMRES will fail to converge. In contrast, nlTGCR(m) and AA do not exclusively rely on one
Jacobian at each iteration.

 \subsection{Image classification using ResNet}\label{sec:resnet}
\begin{figure}[htbp]
   \centering
   \subfigure[Training Loss]{
     \includegraphics[width=0.45\textwidth]{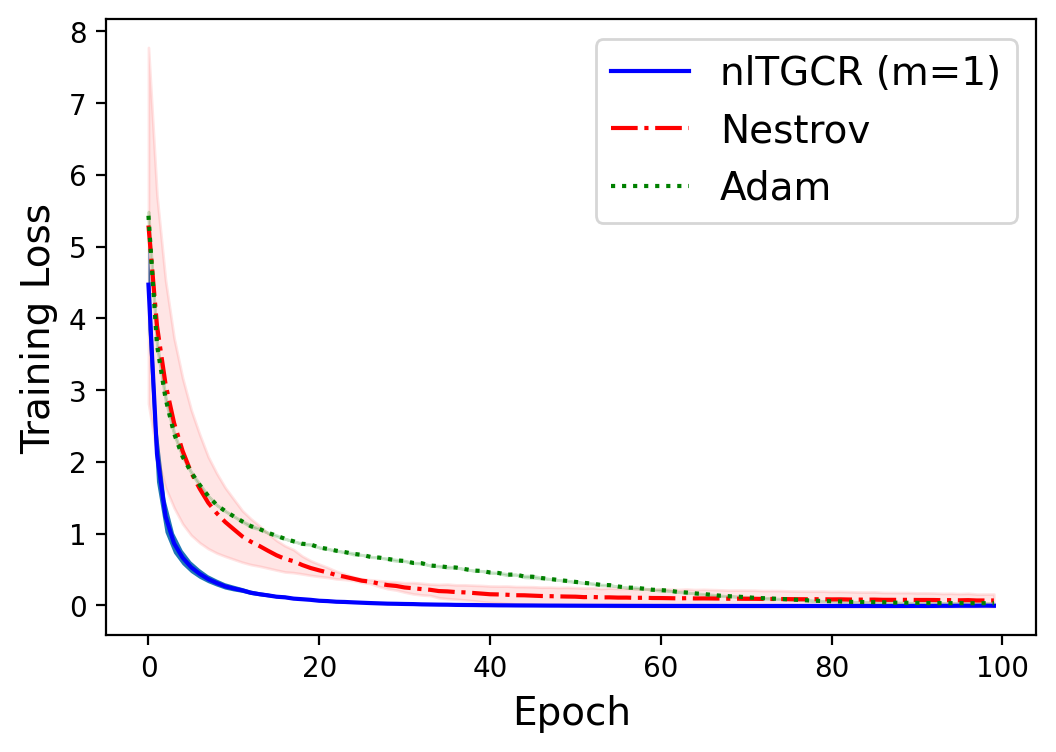}
     \label{fig:ra}
   }
   \subfigure[Test Accuracy]{
     \includegraphics[width=0.45\textwidth]{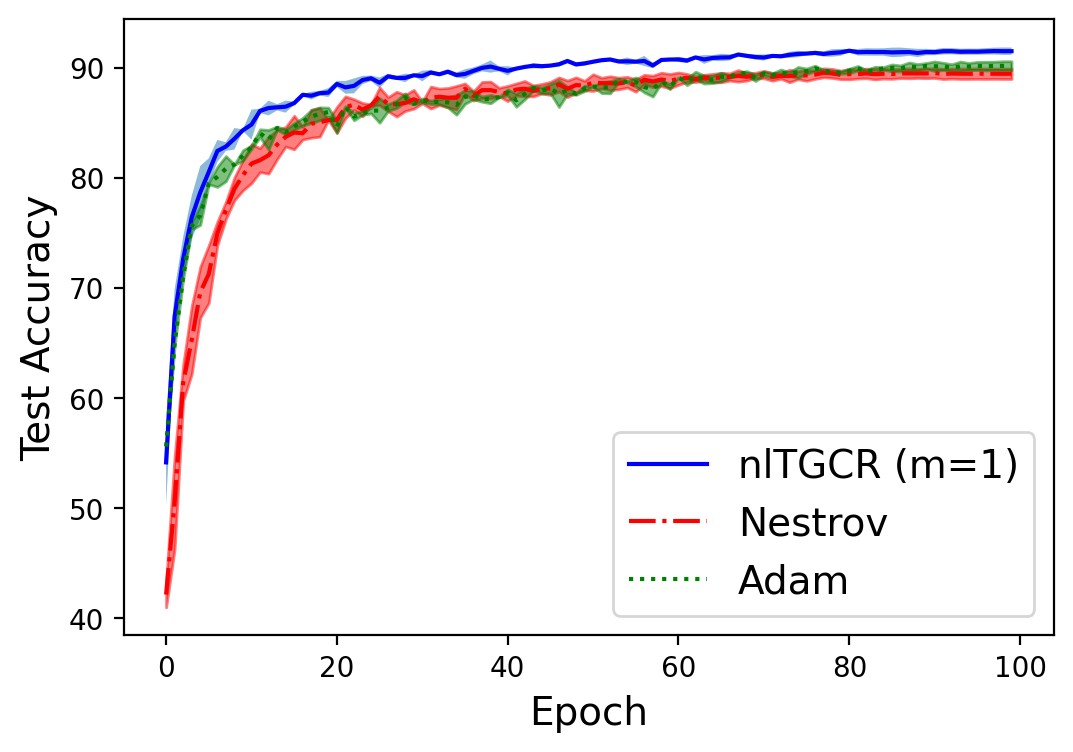}
     \label{fig:rb}
   }
   \caption{Image Classification on CIFAR10 using ResNet (Averaged over 5 independent runs). nlTGCR(m=1), Adam, and momentum achieved a test accuracy of $91.56\%, 90.13 \%, 89.53 \%$ respectively. }
   \label{fig:resnetfig}
	\end{figure}

 We now test the usefulness of nlTGCR(m) for deep learning applications by first
 comparing it with two commonly used optimization algorithms, Adam
 \cite{kingma2017adam} and momentum.  We report the training mean squared error (MSE) and test
 accuracy on the CIFAR10 dataset \cite{Krizhevsky2009LearningML} using ResNet
 \cite{he2015deep}. We employed a ResNet 18 architecture from PyTorch 
 \footnote{https://github.com/pytorch/vision/blob/main/torchvision/models/resnet.py}.
 The window size $m$ is 1 for nlTGCR(m) because we found a large window size did
 not help improve the convergence too much in our preliminary experiments. The
 hyperparameters of the baseline methods are set to be the same as suggested in
 \cite{he2015deep}, i.e., the learning rate is 0.001 and 0.1 for Adam and
 momentum, respectively.  Figure \ref{fig:ra} depicts the training loss over
 time for each optimization algorithm, and Figure \ref{fig:rb} shows the
 corresponding test accuracy. As can be seen nlTGCR($m=1$) achieved the best
 performance in terms of both convergence speed and accuracy. The experimental
 results revealed that nlTGCR outperformed the baseline methods by
 a significant margin. It is worth noting that nlTGCR converges to a loss of 0
 for this problem, which empirically verifies the theoretically established global convergence  properties of
 the method.  This suggests that nlTGCR($m=1$) may lead to an interesting
 alternative to the state-of-the-art optimization methods in deep learning.

\subsection{Learning dynamic using Neural-ODE solver}
We conducted experiments using a Neural-ODE solver \cite{chen2019neural}
to learn underlying dynamic ODEs from sampled data. In our work, we focused on studying the spiral
function {\[
\frac{dz}{dt} = 
\begin{bmatrix}
-0.1 & -1.0 \\
1.0 & -0.1
\end{bmatrix}
z,
\]}which is a
challenging dynamic to fit and is often used as a benchmark for testing
the effectiveness of machine learning algorithms. We generated the
training data by randomly sampling points from the spiral and adding
small amounts of Gaussian noise. The goal was to train the model to
generate data-like trajectories.

\begin{figure}[htbp]
		\centering
		\subfigure[Training Loss]{
			\includegraphics[width=0.3\textwidth]{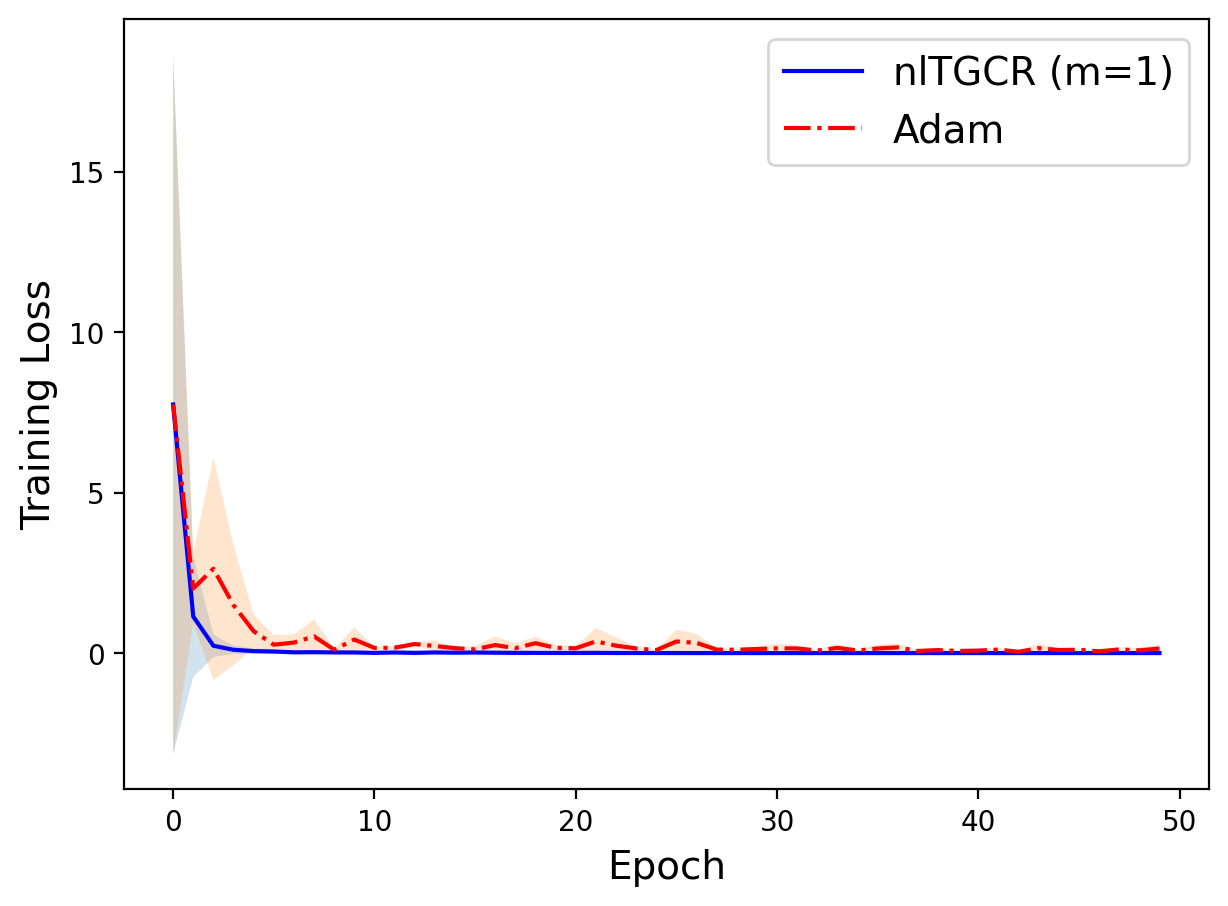}
			\label{fig:na}
		}
		\subfigure[Adam after 50 epochs]{
		\includegraphics[width=0.3\textwidth]{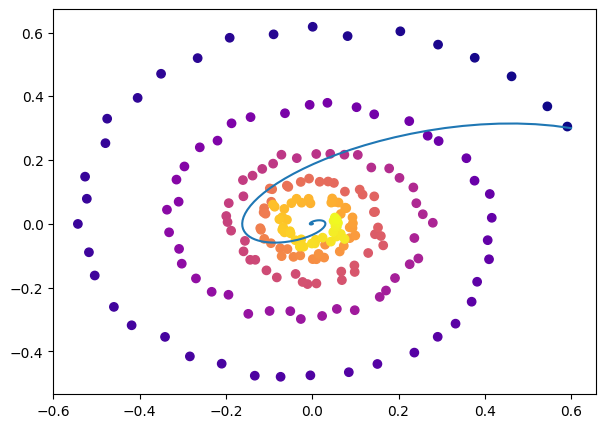}
		\label{fig:nb}
	}
		\subfigure[nlTGCR(m=1) after 15 epochs]{
			\includegraphics[width=0.3\textwidth]{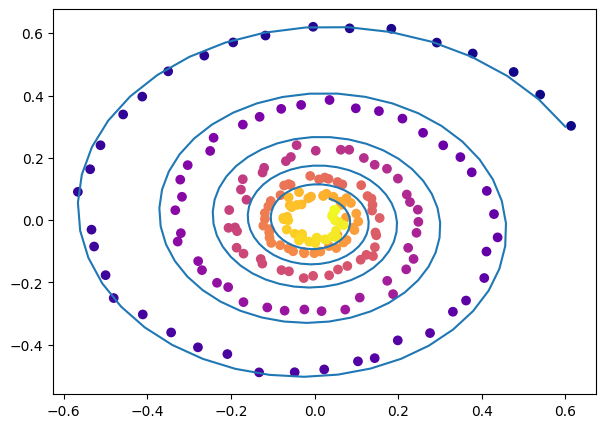}
			\label{fig:nc}
		}
		\caption{Learning true dynamic using Neural ODE. It can be observed that nlTGCR(m=1) converges much faster and can recover the true dynamic up to 5$\times$ faster than Adam.}
		\label{fig:NodeFig}
	\end{figure}

However, training such a model is computationally expensive. Therefore, we
introduced nlTGCR(m) to recover the spiral function quickly and accurately, with
the potential for better generalization to other functions. Our experiments
involved training and evaluating a neural-ODE model on the sampled dataset for
the spiral function compared with Adam and momentum. After a grid-search, we set
the learning rate as 0.01 for Adam and window size $m=1$ for nlTGCR(m). We reported
the MSE training loss and visualized the model's ability to
recover the dynamic in Figure \ref{fig:NodeFig}. We did not visualize the {momentum results} as it took over 100 epochs to converge.

Figure \ref{fig:na} shows that
nlTGCR(m=1) converges faster and more stably than Adam. Additionally, Figures
\ref{fig:nb} and \ref{fig:nc} demonstrate that nlTGCR(m=1) is capable of generating
data-like trajectories in just 15 epochs, whereas Adam struggles to converge
even after 50 epochs. This experiment demonstrates the superiority of nlTGCR($m=1$)
for this interesting application,  relative to commonly used optimizers
such as Adam and momentum. 
	
	\subsection{Node classification using Graph Convolutional Networks}
	\label{sec:gcn}
	Finally, we explore the effectiveness of nlTGCR(m) in deep learning
        applications by implementing graph convolutional networks
        (GCNs)~\cite{kipf2016semi}. We use the commonly used Cora
        dataset~\cite{mccallum2000automating} which contains 2708 scientific
        publications on one of 7 topics. Each publication is described by a
        binary vector of 1433 unique words indicating the absence or presence in
        the dictionary. This network consists of 5429 links representing the
        citation. The objective is a node classification via words and
        links. The neural network has one GCN layer and one dropout layer of
        rate 0.5. We adopted the GCN implementation directly from 
        PyTorch-Geometric~\cite{pyg}.

\begin{figure}[htbp]
	\centering
	\subfigure[Training Loss]{
		\includegraphics[width=0.45\textwidth]{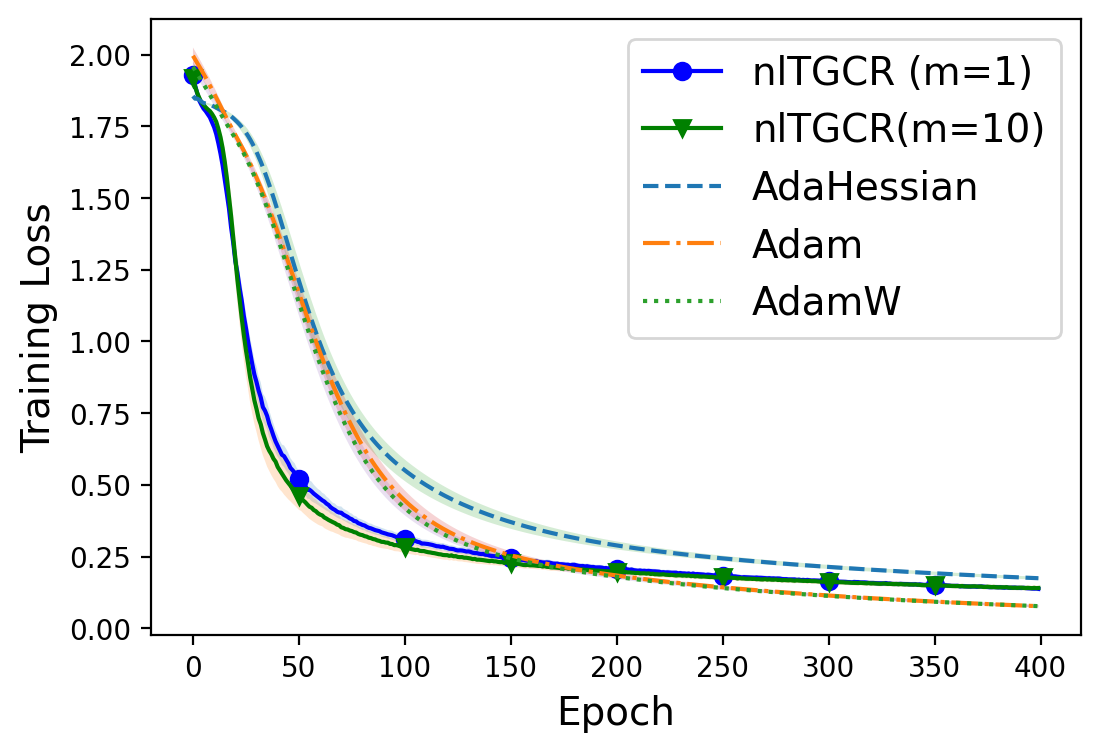}
		\label{fig:2a}
	}
	\subfigure[Test Accuracy]{
		\includegraphics[width=0.45\textwidth]{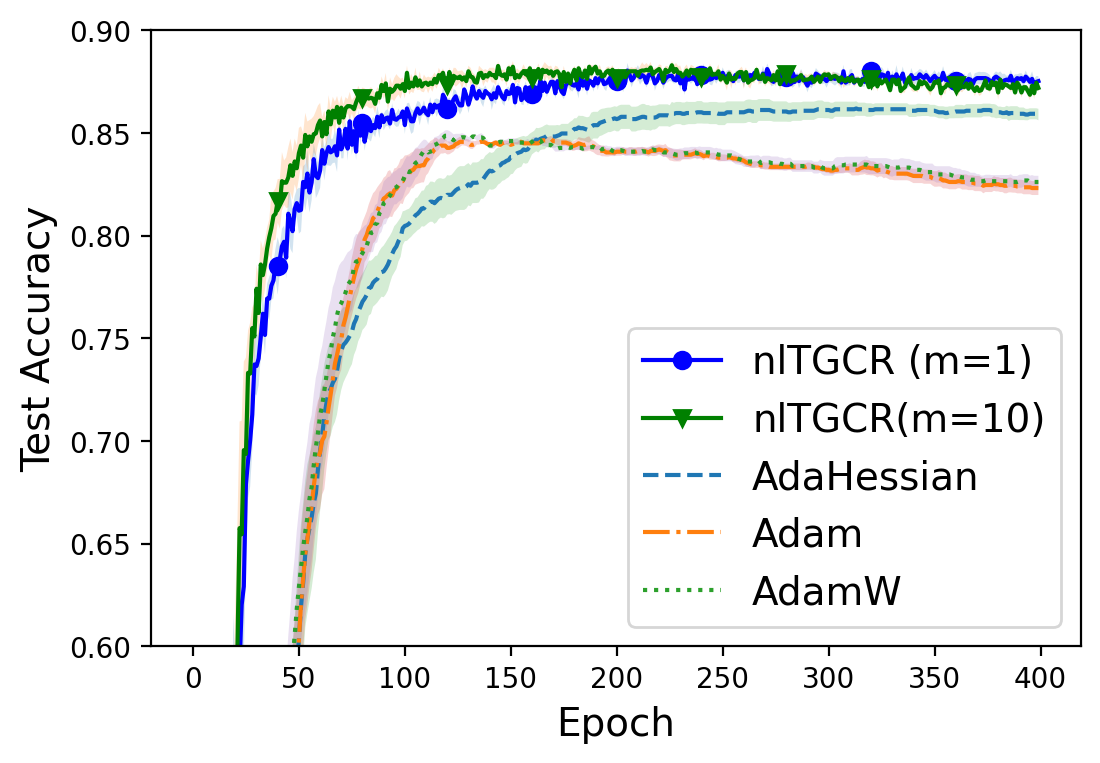}
		\label{fig:2b}
	}
	\caption{GCN on Cora. nlTGCR(m=1) achieved a test accuracy of $88.13\%$, which is $4.02\%$ higher than the second best baseline AdamW.}
	\label{fig:Corafig}
\end{figure}

We set {$m=1$ and $m=10$} in nlTGCR(m) and  compare it with Adam~\cite{kingma2014adam}, AdamW
\cite{loshchilov2019decoupled}, and AdaHessian \cite{yao2020adahessian} with
learning rates $lr=0.01, 0.01, 0.003$ respectively after a grid search. We
present results of the training loss (Cross Entropy) and test accuracy in Figure
\ref{fig:2a} and \ref{fig:2b}. Note that a lower loss function does not
necessarily means a better solution or convergence because of the bias and
variance trade-off. Although Adam achieved a lower training loss, it is not
considered as a better solution because its poor generalization capability on
unseen datasets.  We can observe that nlTGCR($m=1$) shows the best performance in
this task. Specifically, it is not necessary to use a large window size since
there is no significant difference between $m=1$ and $m=10$.  This experiment
shows again that nlTGCR($m=1$) can be adapted to the solution of {deep learning} problems.

\section{Concluding remarks} 
The initial goal of this study was primarily to seek to develop Anderson-like
methods that can take advantage of the symmetry of the Hessian in optimization
problems. What we hope to have conveyed to the reader is that, by a careful
extension of linear iterative schemes, one can develop a whole class of methods,
of which nlTGCR(m) is but one member, that can be quite effective, possibly more so
than many of the existing techniques in some situations.  We are cautiously
encouraged by the results obtained for deep learning problems although much work
remains to be done to adapt and further test nlTGCR(m) for the stochastic context.
In fact, our immediate research plan is precisely to perform such an in-depth
study that focuses on {deep learning} applications.

\section*{Acknowledgement} 
We would like to thank the editor and the two reviewers for their insightful suggestions that have significantly enhanced the quality of our paper.

\bibliographystyle{siam}

\bibliography{strings,local,local1,saad,accelbib,accel1,KRYLOV,msecant}

\newpage

\section{Appendix~A: More on GCR for solving linear systems}
\label{sec:Appendix}
\newcommand{\LL}{\mathcal{L}} A number of results on the GCR algorithm
for linear systems are known but their statements or proofs are not
readily available from a single source. For example, it is intuitive,
and well-known, that when $A$ is Hermitian then GCR will simplify to
its CG-like algorithm known as the Conjugate Residual algorithm, but a proof
cannot be easily found.  This section summarizes the most important
ones of these results with an emphasis on a unified presentation that
exploits a matrix formalism.

\subsection{Conjugacy and orthogonality relations} 
We  start with Algorithm \ref{alg:gcr} where we assume no truncation
($m  = \infty $). It is convenient to use a matrix formalism
for the purpose of unraveling some relations and so we start by defining:
\eq{eq:Rmat}
R_k = [r_0, r_1, \cdots, r_k] .
\en 
The relation in Line 7 of the algorithm
can be rewritten in matrix form as follows:
\eq{eq:HesRel}
R_{k} = P_{k}  B\up{k},
\en
where $ B\up{k}$ is an upper triangular matrix of size $(k+1) \times (k+1)$ defined as follows
\eq{eq:bkp1}
B\up{k} \in \RR^{(k+1)\times (k+1)}, \qquad 
B\up{k}_{ij} = \left\{
  \begin{array}{lcl}
    0 & \text{if} & i>j \\
    1 & \text{if} & i=j \\
    \beta_{(i-1),(j-1)}    &\text{if} & i<j
  \end{array} \right.
\en
Similarly, note that the relation from Line~6 of the
algorithm can be recast as:
\eq{eq:AP2R} AP_k = R_{k+1}  \underline H\up{k},
\en
where $\underline H\up{k}$ is a $(k+2)\times (k+1) $ lower bidiagonal
matrix with 
\eq{eq:uHk}
 \underline H\up{k}_{ij} = \left\{
  \begin{array}{lcl}
    \frac{1}{\alpha_{j-1}} & \text{if} & i=j \  \\
    \frac{-1}{\alpha_{j-1}} & \text{if} & i=j+1 \\
            0           & \text{if} & i<j \ \text{or} \ i>j+1 
  \end{array} \right.
\en

\begin{proposition} (Eisenstat-Elman-Schultz \cite{Eis-Elm-Sch})
  \label{prop:semiconj}
The residual vectors produced by (full) GCR are semi-conjugate.
\end{proposition}

\begin{proof}
  Semi-conjugacy means  that
  each $r_j$ is orthogonal to $Ar_0, A r_1, \cdots Ar_{j-1}$,
  so we need to show that  $ R_k^T A R_k$ is upper triangular.
  We know that each $r_{j+1}$ is orthogonal to
$Ap_i$ for $i=1,\cdots , j$ a relation we write as
\eq{eq:RTAP}
R_{k}^T A P_{k} =  U_k
\en
where $U_k$ is some upper triangular matrix.
Then, from \eqref{eq:HesRel} we have
$  A R_k = A P_k B\up{k} $ and therefore:
\eq{eq:RTAR}
   R_k^T A R_k = R_k^T A P_k B\up{k}  =   U_k B\up{k} 
\en
which is upper triangular as desired.
\end{proof} 

We get an immediate consequence of this result for the case when $A$ is
Hermitian.

\begin{corollary}
   When $A$ is symmetric real (Hermitian complex) then
   the directions $\{r_j\}$, are $A$-conjugate.
 \end{corollary}
\begin{proof}
Indeed, when $A$ is symmetric real the matrix
$R_k^T A R_k =   U_k B\up{k} $ is also symmetric and since it is
upper triangular it must be diagonal which shows that the $r_j$'s are
$A$-conjugate.
\end{proof}

In this situation, we can write
\eq{eq:Rconj}
R_k^T A R_k = D_k ,
\en
where $D_k$ is a $(k+1)\times (k+1)$ diagonal matrix. The next result shows that the algorithm
simplifies when $A$ is symmetric. Specifically, the scalars $\beta_{ij}$ needed in the
orthgonalization in Line 7 are all zero except $\beta_{jj}$.

\begin{proposition}When $A$ is symmetric real, then the matrix
  $(AR_k)^T (AP_k)$ is lower bidiagonal.
  \label{prop:cr}
\end{proposition} 
\begin{proof} 
As a result of \nref{eq:AP2R}  the matrix $(AR_k)^T (AP_k)$ is:
\[
(AR_k)^T (AP_k) = (AR_k)^T  R_{k+1} \underline H\up{k} . 
\]  
Observe that $(AR_k)^T  R_{k+1} = [D_k, 0_{(k+1)\times 1}] $ and the product
$ [D_k, 0_{(k+1)\times 1}] \underline H\up{k}$ yields the $(k+1)\times (k+1)$  matrix obtained
from $ \underline H\up{k} $ but deleting its last row which we denote by
$H\up{k}$. Therefore,
\[
(AR_k)^T (AP_k) = D_k H\up{k} 
\]  
is a $(k+1)\times (k+1)$ bidiagonal matrix.
\end{proof} 

\subsection{Break-down of GCR} 
Next we examine the conditions under which the full GCR breaks down.
\begin{proposition}
  When $A$ is nonsingular, the only way in which (full) GCR
  breaks down  is when it produces an exact solution. In other words its only
  possible breakdown is the  so-called `lucky breakdown'.
\end{proposition}
\begin{proof}
  The algorithm breaks down only 
  if $Ap_{j+1}$ produced in Line~7 is zero, i.e., when $p_{j+1} == 0$ since $A$
  is nonsingular.
  In this situation $r_{j+1} = \sum_{i=0}^j \beta_{ij} p_i$.
  It can be easily seen that each $p_i $ is of the form
  $p_i = \mu_i (A) r_0$ where $\mu_i$ is a polynomial of degree $i$.
  Similarly $r_{j+1} = \rho_{j+1} (A) r_0$ in which 
  $\rho_{j+1}$ is a polynomial of degree (exactly) $j+1$.
  Therefore, the polynomial
  $\pi_{j+1}(t) \equiv \rho_{j+1}(t)-\sum_{i=0}^j \beta_{ij} \mu_i(t)$
  is a  polynomial of degree exactly $j+1$ such that
  $\pi_{j+1}(A)r_0=0 $. Thus, the degree of the  \emph{minimal polynomial} for
  $r_0$ is $j+1$ and we are in the standard situation of a lucky breakdown.
 Indeed, since the algorithm produces a solution
  $x_{j+1}$ that  minimizes   the residual norm, and because
  the degree of the minimal polynomial for $r_0$ is equal to $j+1$
  we must have   $r_{j+1} = 0$. 
\end{proof}

Note that the proposition does not state anything about convergence. The
iterates that are computed will have a residual norm that is non-increasing
but the iterates may stagnate.  Convergence can be shown in the case when $A$ is
positive definite, i.e., when its symmetric part is SPD.

In addition, the proof of this result requires that the solution that is
produced has a minimal residual which is not the case for the truncated version.
Thus, in the truncated version,
we may well have a situation where $ p_{j+1} = 0$ but the 
solution $x_{j+1}$ is not exact. If we had $ p_{j+1} = 0$ it would mean that
$\eta_{j+1}(A)r_0 = 0$, i.e., the minimal polynomial for $r_0$ is again of
degree exactly $j+1$. This is an unlikely event that we may term `unlucky
breakdown'. However, in practice, the more common situation that can take place
is to get a vector $p_{j+1}$ with a small norm.


Suppose now that $A$ is positive definite, i.e. that its symmetric part is SPD.
Let us assume that $p_k = 0$ but $r_k \ne 0$ -- which represents the scenario of
an `unlucky breakdown' at step $k$.
Then since the last column of $P_k$ is a zero vector 
the last column of the matrix $U_k$ in \eqref{eq:RTAP} is also zero.
This in turn would imply that the last \emph{row} of the product $  U_k B\up{k}$
in \eqref{eq:RTAR} is  zero. 
However, this can't happen because according to \eqref{eq:RTAR} 
it  is equal to the last row of the matrix
$R_k ^T A R_k$ where $A$ is positive definite.
Thus, the `unlucky breakdown' scenario invoked above is only possible when
$A$ is indefinite.

\subsection{Properties of the induced approximate inverse}
When  $x_0 = 0$ the approximate solution obtained at the end of the
algorithm is of the form $x_{k+1} = P_k V_k^T b$.  We say that the algorithm induces
the approximation $B_k = P_k V_k^T$ to the inverse of $A$. Note that even in the
case when $A$ is symmetric, $B_k$ is not symmetric in general. However, $B_k$
obeys a few easy-to-prove properties stated next. 
\begin{proposition}
  Let $\LL_k = \Span (V_k)$ and let 
   $\pi = V_k V_k^T$ be the orthogonal projector onto $\LL_k$.
  The (full) GCR algorithm induces the approximate inverse $B_k = P_k V_k^T$ which
  satisfies the following properties:
  \begin{enumerate}
    \item $B_k = A\inv \pi $
  \item $B_k $ inverts $A$ exactly in $\LL_k$, i.e.,    $B_k x = A \inv x$
    for $ x \in \LL_k$. Equivalently, $B_k \pi = A\inv \pi$.
  \item $A B_k$ equals the orthogonal projector $ \pi$.
  \item When $A$ is symmetric then
    $B_k$ is self-adjoint when restricted to $\LL_k$. 
  \item $B_k A x = x$ for any $x  \in \Span \{P_k \}$, i.e., 
    $B_k$ inverts $A$ exactly from the left when $A$ is restricted to
    the range  of $P_k$.

  \item $B_k A$ is the projector onto $\Span \{P_k \}$ and orthogonally to $A^T \LL_k$.
  \end{enumerate}
\end{proposition}
\begin{proof}
  (1)
  The first property follows from the relation $A P_k = V_k$ and the definition of $B_k$.

(2) To prove the second property we write a member of $\Span (V_k)$ as $x = V_k y$. Then from
the previous property we have
\[B_k x = A\inv V_k V_k^T V_k y = A\inv V_k y =  A\inv x . \]

(3) $A B_k = A P_k V_k^T = V_k V_k^T = \pi$.

(4) The self-adjointness of $B_k$ in $\LL_k$ is a consequence of (2). It can also be
readily verified as follows. For any vectors $x, y \ \in \ \LL_k$ write:
\[
  (B_k x, y) = (A \inv \pi x, y) = (A \inv x, y) = (x, A \inv y) = (x, A \inv \pi y) =
  (x, B_k y) .
\]


(5) Let a member of $\Span (P_k) $ be written as $x = P_k y$ and apply $B_k A $ to $x$:
\[
  B_k A x = B_k A P_k y = B_k V_k y = A\inv V_k V_k^T V_k y =  A\inv V_k y =  P_k y = x . \]
Therefore, $B_k A$ leaves vectors of $\Span (P_k)$ unchanged.

(6) Because $B_k A$ leaves vectors of $\Span (P_k)$ unchanged it is a projector, call it $\pi_O$ (for oblique),
with $\text{Ran} (\pi_O) = \Span \{P_k \}$.
We now need to show that $ (u - \pi_O u) \perp A^T \LL_k $ for any $u$. Since $A^T V_k$ is a basis for
$A^T \LL_k$, this is  equivalent to
showing that  $  (A^T V_k)^T (u - \pi_O u)) = 0$ for any $u$. We have for any given vector $u \in \RR^n$
\begin{align*}
(A^T V_k)^T (u - \pi_O u) & = V_k^T (A u - A P_k V_k^T A u) = V_k^T (A u  -V_k V_k^T A u) \\
& = V_k^T (I  -V_k V_k^T)  A u  = 0 . 
\end{align*}
\end{proof}

Thus (4) and (6) indicate that while $A B_k$ is an orthogonal projector, $B_k A$ is an oblique projector.
  Although (4) is an obvious consequence of (3), it is interesting to note that it
  is rather similar to relations obtained in the context of Moore-Penrose   pseudo-inverses.


  \section{Appendix~B: Convergence analysis}\label{sec:alterproof}
  We provide an alternative analysis of the global convergence of nlTGCR based on the backtracking line search
  strategy. We will make the  following assumptions:

\bigskip\noindent\textbf{Assumption~A:}
\eq{asa}
\text{ The  set} ~ S = \{x: \phi(x) \leq \phi(x_0)\} ~\text{is bounded};
\en
\bigskip\noindent\textbf{Assumption~B:}
\eq{asb}
\nabla \phi  ~\text{is L-Lipschitz},\|\nabla \phi(x) - \nabla \phi(y)\| \leq L\|x-y\|, ~0<L<\infty.
\en
\bigskip\noindent\textbf{Assumption~C:} There exists a $\gamma>0$ such that 
$d_j$ produced by Algorithm \ref{alg:nltgcr} satisfies
\eq{asc}
-\frac{d_j^{\top}\nabla \phi(x_j)}{\|d_j\|\|\nabla\phi(x_j)\|}\geq \gamma >0,  \forall j
\en
\bigskip\noindent\textbf{Assumption~D:} There exist two positive constants
$\mu$ and $\Theta $ such that,
\eq{asd}
\|d_j\|\geq \mu \|\nabla\phi(x_j)\|, \quad \|d_j\|\leq \Theta, \forall j.
\en
\begin{theorem}
	For any scalar function $\phi(x)$, consider the iterates $x_{j+1} = x_{j} + \alpha_{j}d_{j}$ with descent directions $d_j$ produced by Algorithm \ref{alg:nltgcr} and stepsizes $\alpha_j$ produced by backtracking line search \eqref{eq:Armijo3}. Suppose assumptions (A-D) also hold, then 
	\eq{convergence to stationary point}
	\lim_{j\rightarrow\infty} \|\nabla \phi(x_j)\| = 0.
	\en
\end{theorem}
This theorem  is a well-known result in optimization and it
can be found in various sources, including \cite{griva2009linear}. We omit its proof. 

\end{document}